\documentclass{amsart}
\usepackage{amssymb}
\newcommand{\card}{\mathrm{Card}}

\newcommand{\Q}{{\mathbb Q}}
\newcommand{\Z}{{\mathbb Z}}
\newcommand{\bN}{{\mathbb N}}
\newcommand{\N}{\mathcal N}

\newcommand{\R}{{\mathbb R}}
\newcommand{\F}{{\mathbb F}}

\newcommand{\OO}{{\mathcal O}}
\newcommand{\Pp}{{\mathfrak{P}}}

\newcommand{\ga}{{\mathfrak a}}
\newcommand{\gb}{{\mathfrak b}}

\newcommand{\gd}{{\mathfrak d}}
\renewcommand{\gg}{{\mathfrak g}}
\newcommand{\gp}{{\mathfrak p}}
\newcommand{\gq}{{\mathfrak q}}
\newcommand{\gl}{{\mathfrak l}}

\newcommand{\Tc}{{\mathcal T}}

\newcommand{\eps}{{\varepsilon}}
\newcommand{\M}{{\rm M}}
\newcommand{\h}{{\rm h}}
\newcommand{\cA}{{\mathcal A}}

\newcommand{\cL}{{\mathcal L}}
\newcommand{\cD}{{\mathcal D}}
\newcommand{\bC}{{\mathbb C}}
\newcommand{\q}{{\chi_0}}

\newcommand{\kro}[2]{\left( \frac{#1}{#2} \right) }

\begin {document}
\newtheorem{thm}{Theorem}
\newtheorem{lem}{Lemma}[section]

\newtheorem{cor}[lem]{Corollary}

\newtheorem{prop}[lem]{Proposition}

\newtheorem*{defn}{Definition}

\theoremstyle{definition}

\newtheorem{ex}{Example}

\theoremstyle{remark}
\newtheorem*{remark}{Remark}

\title[Classical and Modular Approaches]{Classical
and Modular Approaches to Exponential Diophantine Equations \\
II. The Lebesgue--Nagell Equation} 
\author{Yann Bugeaud, Maurice Mignotte, Samir Siksek}
\address{Yann Bugeaud\\
	Universit\'{e} Louis Pasteur\\
	U. F. R. de math\'{e}matiques\\
	7, rue Ren\'{e} Descartes\\
	67084 Strasbourg Cedex\\
	France}
\email{bugeaud@math.u-strasbg.fr}

\address{Maurice Mignotte\\
        Universit\'{e} Louis Pasteur\\
        U. F. R. de math\'{e}matiques\\
        7, rue Ren\'{e} Descartes\\
        67084 Strasbourg Cedex\\
        France}
\email{mignotte@math.u-strasbg.fr}

\address{Samir Siksek\\
	Department of Mathematics and Statistics\\
        College  of Science \\
        Sultan Qaboos University \\
        P.O. Box 36 \\
        Al-Khod 123 \\
        Oman}
\email{siksek@squ.edu.om}

% \dedicatory{}
\date{\today}
\thanks{S. Siksek's work is funded by a grant from Sultan Qaboos University (IG/SCI/DOMS/02/06)}

\keywords{Diophantine equations, Ramanujan--Nagell, Frey curves, level-lowering, linear forms in logarithms, Thue equations}
\subjclass[2000]{Primary 11D61, 11J86, Secondary 11D59, 11Y50}

\begin {abstract}
This is the second in a series of papers where we combine the
classical approach to exponential Diophantine equations (linear
forms in logarithms, Thue equations, etc.) with a modular 
approach based on some of the ideas of the proof of Fermat's Last 
Theorem. In this paper we give a general and powerful lower bound
for linear forms in three logarithms. We use this lower bound,
together with a combination of classical, elementary and substantially
improved modular methods to solve completely the Lebesgue-Nagell equation
\[
	x^2+D=y^n, \qquad \text{$x$, $y$ integers, $n\geq 3$}, 
\]
for $D$ in the range $1 \leq D \leq 100$.
\end {abstract}
\maketitle

\section{Introduction}
Arguably, the two most celebrated achievements of the 20th century
in the field of Diophantine equations have been Baker's theory
of linear forms in logarithms, and Wiles' proof of Fermat's Last Theorem.
We call Baker's approach to Diophantine equations the \lq classical approach\rq.
The proof of Fermat's Last Theorem is based on what we term the \lq modular approach\rq.
The proponents of the classical approach are too many to mention; the modular
approach is still in its infancy, but among the early contributers let us
just mention Frey, Serre, Ribet, Darmon, Merel, Kraus, Bennett, Skinner, 
Ivorra, etc.  %$%

The motivation for our series of papers, of which this is the second,
is that neither approach (on its own, and as it stands at the moment)
is powerful enough to resolve unconditionally many of the outstanding
exponential Diophantine equations. Our thesis is that one should, %#%
where possible,  attack 
exponential Diophantine equations by a combination of the 
classical and modular approaches. 
The precise aims of this series were formulated in
our first paper~\cite{BMS} as follows:
\begin{enumerate}
\item[(I)] To present theoretical improvements to various aspects of
the classical approach.
\item[(II)] To show how local information obtained through the modular
approach can be used to reduce the size of the bounds, both for
exponents and for variables, of solutions to exponential Diophantine
equations.
\item[(III)] To show how local information obtained through the modular
approach can be pieced together to provide a proof that there are no
missing solutions less than the bounds obtained in (I), (II).
\item[(IV)] To solve various famous and hitherto outstanding exponential
Diophantine equations.
\end{enumerate}

In \cite{BMS} we gave a lower bound for linear forms in three logarithms,
and used a combination of classical and modular methods to determine
all the perfect powers in the Fibonacci and Lucas sequences.
In this paper, we give a new lower bound for linear forms in three logarithms
that is more general and powerful than the one given in the previous paper. 
We are also concerned with the following exponential Diophantine equation,
which we call the Lebesgue--Nagell equation:
\begin{equation}\label{eqn:first}
	x^2+D=y^n, \qquad \text{$x,~y$ integers, $n \geq 3$}.
\end{equation}
Here, $D$ denotes a non-zero integer. %$%
The reason for the name Lebesgue--Nagell is given in Section~\ref{sec:survey}, 
together with some historical remarks.
But for now we mention
that the equation had previously been solved for $81$ values of $D$
in the range $1 \leq D \leq 100$, using elementary, classical and modular methods;
the remaining values are clearly beyond these methods as they stand.
In this paper we apply our lower bound for linear forms in three logarithms, 
together with a combination of elementary, classical,
and substantially improved modular methods to prove the following Theorem.
\begin{thm}\label{thm:main}
All solutions to equation~(\ref{eqn:first}) with $D$ in the range
\begin{equation}\label{eqn:Drange}
	1 \leq D \leq 100
\end{equation}
are given in the Tables at the end. In particular, the only %$%
integer solutions $(x, y, n)$ to the equation
\[
x^2+7=y^n, \qquad n \geq 3,
\]
satisfy $|x|=1,~3,~5,~11,~181$.
\end{thm}

We choose to give a complete proof of Theorem 1, rather than treating %$%
the 19 remaning values of $D$ in the range~(\ref{eqn:Drange}) .

It is noted that the solutions for even $n$ can be 
deduced quickly, for then $D$ is expressible
as a difference of squares. It is therefore sufficient to solve the 
equation
\begin{equation}\label{eqn:p}
        x^2+D=y^p, \qquad \text{$x,~y$ integers, $p \geq 3$ is prime};
\end{equation}
the solutions to~(\ref{eqn:first}) can then be recovered from the 
solutions to~(\ref{eqn:p}).

We give three modular methods for attacking~(\ref{eqn:p}). Two are refinements
of known methods, and a third that is completely new. Using a computer
program based on these modular methods, we can show -- for 
any $D$ in the above range -- that the exponent $p$ is large (showing that
$p > 10^9$ is quite practical). Our modular approach 
also yields the following rather surprising result: 
either each prime factor of $y$ divides $2D$, or $y > (\sqrt{p}-1)^2$.
We are then able to deduce not only that $p$ is large, but also that $y$ is
large. This information helps to reduce the size of the upper bound on $p$ obtained 
from the lower bound for the linear forms in three logarithms, making 
the computation much more practical.  Our total computer time for the 
computations in this paper is roughly $206$ days on various workstations
(the precise details are given in due course).

Using our approach should make it possible to solve~(\ref{eqn:first})
for any reasonable $D$ that is {\bf not} of the form $D=-a^2 \pm 1$;
if $D$ is of this form then the equation~(\ref{eqn:first})
has a solution $(x,y)=(a,\pm 1)$ for all odd values of the exponent $n$,
and the modular methods we explain later are not very successful in
this situation. To deal with this case requires further considerations
which we leave for another paper. 
Notice, however, that we solve the case $D=1$.

%This is the reason why the value $D=1$
%is omitted from our range~(\ref{eqn:Drange}). 

We would like to warmly thank Mihai Cipu for pointing our many imperfections
in a previous version of this paper, and Guillaume Hanrot for help with
solving Thue equations.

\section{On the History of the Lebesgue--Nagell Equation}\label{sec:survey}

Equation~(\ref{eqn:first}) has a long and glorious history, and there
are literally hundreds (if not thousands) of papers devoted to special
cases of this equation. Most of these are concerned with
equation~(\ref{eqn:first}) either for special values of $n$ 
or special values of $y$.
For example, for $D=2,~n=3$, Fermat asserted that he
had shown that the only solutions are given by $x=5,~y=3$;
a proof was given by Euler \cite{Euler}. Equation~(\ref{eqn:first})
with $n=3$ is the intensively studied Mordell equation (see~\cite{GPZ} 
for a modern approach).

Another notable special case is the generalized Ramanujan--Nagell
equation
\begin{equation}\label{eqn:gRN}
x^2+D=k^n,
\end{equation}
where $D$ and $k$ are given integers.
This is an extension of the Ramanujan--Nagell equation $x^2+7=2^n$,
proposed by Ramanujan~\cite{Ramanujan} in 1913 and
first solved by Nagell~\cite{Nagell} in 1948 (see also the collected papers
of Nagell \cite{Nagell2}).
This equation has exactly five solutions with $x \ge 1$ and is in this
respect singular: indeed, Bugeaud and Shorey~\cite{BuSho}
established that equation~(\ref{eqn:gRN}) with $D$ positive
and $k$ a prime number not dividing $D$
has at most two solutions in positive integers $x,~n$,
except for $(D, k) = (7, 2)$. They also list all the pairs
$(D, k)$ as above for which equation~(\ref{eqn:gRN}) has exactly two
solutions. We direct the reader to~\cite{BuSho}
for further results and references.

Returning to equation~(\ref{eqn:first}), the first result 
for general $y$, $n$ seems 
to be the proof in 1850 by V. A. Lebesgue \cite{Lebesgue}
that there are no non-trivial solutions for $D=1$. 
The next cases to be solved
were $D=3,~5$ by Nagell~\cite{Nagell} in 1923. It is for this reason 
that we call equation~(\ref{eqn:first}) the Lebesgue--Nagell equation.
The case with $D=-1$ is particularly noteworthy: 
a solution was sought for many years as a special case of
the Catalan conjecture. This case was finally 
settled by Chao Ko \cite{Chao} in 1965. %$%

The history of the Lebesgue--Nagell equation is 
meticulously documented in an important article by
Cohn~\cite{Cohn1}, and so we are saved the trouble of compiling an
exhaustive survey. In particular, Cohn refines the earlier elementary
approaches of various authors and completes the solution for $77$ values
of $D$ in the range $1 \leq D \leq 100$. The solution for the cases $D=74,~86$
is given by Mignotte and de Weger~\cite{MdW}. %$%
Bennett and Skinner~\cite[Proposition 8.5]{BS}
apply the modular approach to  solve $D=55,~95$. The $19$ remaining
values
\begin{equation}\label{eqn:left}
7, 15, 18, 23, 25, 28, 31, 39, 45, 47, 60, 63, 71, 72, 79, 87, 92, 99, 100,
\end{equation}
 are clearly beyond the scope of Cohn's elementary method, though
Cohn's method can still give non-trivial information even in these cases,
and is revisited in Section~\ref{sec:Cohn}. Moreover, as far as we can see,
the modular method  used by Bennett and Skinner (which is what we
later on call Method I) is not capable of handling these values on it own,
even though it still gives useful information in most cases.

Cohn, in the same paper, also makes a challenge of proving that the only
solutions to the equation
\[
x^2+7=y^n
\]
have $|x|=1,~3,~5,~11,~181$.
This challenge is taken up by Siksek and Cremona \cite{SC}  %$%
who use the modular approach to show that there are no further
solutions for $n \leq 10^8$, nor for composite $n$.
They also suggest that an improvement
to lower bounds in linear forms in three logarithms may
finally settle the problem.
With the benefit of hindsight, we know that they were almost -- though not
entirely -- correct. The 
substantial improvement to lower bounds in linear
forms in three logarithms given here, was certainly needed.
However, for this lower bound to be even more effective, a further insight
obtained from the modular approach was also needed: namely
that $y$ is large as indicated in the introduction.

\section{Reduction to Thue Equations}\label{sec:Thue}

Our main methods for attacking equation~(\ref{eqn:p}) are 
  linear forms in logarithms (to bound the exponent $p$)
and the modular approach, though for some small values
of $p$ it is necessary
to reduce the equation to a family of Thue equations.
The method for reducing equation~(\ref{eqn:p})
to Thue equations is well-known. We do however 
feel compelled to give a succinct recipe for  
this, in order to set up notation that is needed later. 

It is appropriate to point out that there are other
approaches that could be used to solve equation~(\ref{eqn:p})
for small $p$. For $p=3$ we can view the problem as that
of finding integral points on elliptic curve, a problem that
is aptly dealt with in the literature (see \cite{Smart}
and \cite{GPZ}). For $p \geq 5$, the
equation $x^2+D=y^p$ defines a curve of genus $\geq 2$;
one can sometimes determine all rational points on this curve
using the method of Chabauty \cite{CF}, though this would require computing
the Mordell--Weil group of the Jacobian as well 
(see \cite{PS}, \cite{Sch}, \cite{StollI}, \cite{Stoll} and \cite{StollII}). %$%

%We have to deal separately with the cases where $-D$ is a square,
%and $-D$ is a non-square.

%\subsection{The Case $-D$ is a Non-Square}
We do not assume in this section that $D$ is necessarily in
the range~(\ref{eqn:Drange}), merely that $-D$ is not a square.
We write (here and throughout the paper) 
\[
D=D_1^2 D_2, \qquad \text{$D_1,\, D_2$ are integers, $D_2$ square-free}.
\]
Let $\cL=\Q(\sqrt{-D_2})$,
and $\OO$ be its ring of integers.
Throughout the present paper,  %$%
we denote the conjugate of an element $\alpha$ (resp. of an ideal  %$%
$\ga$) by $\overline{\alpha}$ (resp. by $\overline{\ga}$).  %$%
 
Let $\gp_1,\ldots,\gp_r$
be the prime ideals of $\OO$ dividing $2D$. Let $\cA$ be the set of integral
ideals $\ga$ of $\OO$ such that
\begin{itemize}
\item $\ga= \gp_1^{a_1} \cdots \gp_r^{a_r}$, with $0 \leq a_i <p$,
\item $(\ga,\overline{\ga}) \mid 2\, D_1 \sqrt{-D_2}$,
\item the ideal $\ga \overline{\ga}$ is a perfect $p$-th power.
\end{itemize}
If $(x,y)$ is a solution to equation~(\ref{eqn:p}), then
one effortlessly sees that 
\[
(x+D_1 \sqrt{-D_2}) \OO = \ga \gb^p
\]
for some $\ga \in \cA$ and some integral ideal $\gb$.

Now let $\gb_1,\ldots,\gb_h$ be integral ideals forming
a complete set of representatives for the ideal class 
group of $\OO$. Thus $\gb \gb_i$ is a principal ideal
for some $i$, and so  $\gb \gb_i=\beta^\prime \OO$ for some $\beta^\prime \in \OO$.
The fractional ideal $\ga \gb_i^{-p}$ is easily seen to be
also principal. The ideal $\gb$ is unknown, but the
ideals, $\ga,\gb_1,\ldots,\gb_h$ are known. 
We may certainly determine which of the fractional ideals $\ga \gb_i^{-p}$ 
are principal. Let $\Gamma^\prime$ be a set containing
one generator $\gamma^\prime$
for every principal ideal of the form $\ga \gb_i^{-p}$ 
($\ga \in \cA$ and $1 \leq i \leq h$). It is noted that the
elements of $\Gamma^\prime$ are not necessarily integral, but
we know that if $(x,y)$ is a solution to equation
~(\ref{eqn:p}) then 
\[
(x+D_1 \sqrt{-D_2})\OO=\gamma^\prime {\beta^\prime}^p \OO,
\]
for some $\gamma^\prime \in \Gamma^\prime$, and some
$\beta^\prime \in \OO$. 
Finally, define $\Gamma$ as follows: 
\begin{equation*}
	\Gamma=
	\begin{cases}
 \Gamma^\prime  & \text{if $D_2 >0, D_2 \neq 3$, or if $D_2=3$ and $p \neq 3$},\\
 \Gamma^\prime \cup \zeta \Gamma^\prime \cup \zeta^{-1} \Gamma^\prime &
	\text{if $D_2=p=3$, where $\zeta=(1+\sqrt{-3})/2$,}\\
 \cup_{j} \epsilon^j \Gamma^\prime &   %$%
 	\text{if $D_2 <0$, where $j$ ranges over $-(p-1)/2,\ldots,(p-1)/2$,}
	\end{cases}
\end{equation*}
where if $D_2 <0$ (and so $\cL$ is real) we write $\epsilon$
for the fundamental unit.

We quickly deduce the following.
\begin{prop}\label{prop:Thue1}
With notation as above, if $(x,y)$ is a solution to equation~(\ref{eqn:p})
then there exists  $\gamma \in \Gamma$ and $\beta \in \OO$ such that
\[
x+D_1 \sqrt{-D_2} = \gamma \beta^p.
\]
Thus if we let $1$, $\omega$ %%!
be an integral basis for $\OO$ then for some $\gamma \in \Gamma$,
\[
x= \frac{1}{2} \bigl( \gamma (U+V \omega)^p + \overline{\gamma} (U+V\overline{\omega})^p \bigr)
\]
for some integral solution $(U,V)$ to the Thue equation
\[
 \frac{1}{2\sqrt{-D_2}}\bigl( \gamma (U+V \omega)^p - \overline{\gamma} (U+V\overline{\omega})^p \bigr)=D_1.
\]
\end{prop}

\subsection{Results I}
If $q$ is a prime we denote by
$v_q: \Z \rightarrow \Z_{\geq 0} \cup \left\{ \infty \right\}$
the normalized $q$-adic valuation.

We now eliminate all cases where it is inconvenient to carry out level-lowering.
\begin{lem}\label{lem:dmin}
Suppose $D$ is in our range~(\ref{eqn:Drange}). Suppose $(x,y,p)$
is a solution to equation~(\ref{eqn:p}) that is missing from our Tables
at the end. Then $p$ satisfies the following conditions:
\begin{equation}\label{eqn:cond}
\left\{
\begin{array}{ll}
p \geq 7, & \\
p \geq v_q(D)+1 & \text{for all primes $q$,} \\
p \geq v_2(D)+7 & \text{if $v_2(D)$ is even.} \\
\end{array}
\right.
\end{equation}
\end{lem}
\begin{proof}
It is clear that for any particular $D$ there are only a handful of
primes $p$ violating any of these conditions. We wrote a {\tt pari/gp} \cite{pari}
program that solved all the equations~(\ref{eqn:p}) for $p$ violating
~(\ref{eqn:cond}): the program first reduces each such equation
to a family of Thue equations as in Proposition~\ref{prop:Thue1} above.
These are then solved using the in-built {\tt pari/gp} 
function for solving Thue equations 
(this is an implementation of the method of Bilu and Hanrot \cite{BH}). 

It is perhaps worthwhile to record here two
tricks that helped us in this step. First, in writing down the set $\Gamma$
appearing in Proposition~\ref{prop:Thue1} we needed a set of integral
ideals $\gb_1,\ldots,\gb_h$ representing the ideal class group of the quadratic field $\cL$.
Both {\tt pari/gp} and {\tt MAGMA} \cite{magma} have in-built functions that amount to
homomorphisms from the ideal class group as an abstract group,
to the set of fractional ideals, and these can be used to
construct the required set $\gb_1,\ldots,\gb_h$. We have found
however that we get much simpler Thue equations if we search
for the smallest prime ideal representing each non-trivial 
ideal class, and of course taking $1 \OO$ to represent
the trivial ideal class.

To introduce the second trick, we recall that when one
is faced with a Thue equation
\[
a_0 U^p+a_1 U^{p-1} V + \cdots +a_p V^p = b
\]
it is usual to multiply throughout by $a_0^{p-1}$
and make the substitution $U^\prime = a_0 U$,
thus obtaining a monic polynomial on the left-hand side. 
When $a_0$ is large, this greatly complicates the equation.
The second trick is to first search for a unimodular substitution
which makes the leading coefficient $a_0$ small.

After optimizing our program, we were able to complete the
proof in about $22$ minutes on a 1050 MHz UltraSPARC III computer.
\end{proof}

\section{Removing Common Factors}
It is desirable when applying the modular approach to equation~(\ref{eqn:p})
to remove the possible common factors of the three terms in the equation.
This desire leads to a subdivision of cases according to the
possible common factors, as seen in the following elementary Lemma.
Here and elsewhere, for a non-zero integer $a$, the product
of the distinct prime divisors of $a$ is called the radical of $a$, 
and denoted by $\mathrm{rad}(a)$.

\begin{lem}\label{lem:remove}
Suppose that $(x,y,p)$ is a solution to equation~(\ref{eqn:p})
such that $y \neq 0$ and $p$ satisfies the condition~(\ref{eqn:cond}).
Then there are integers $d_1,~d_2$ such that the following conditions
are satisfied:
\begin{enumerate}
\item[(i)] $d_1 > 0$,
\item[(ii)] $D=d_1^2 d_2$,
\item[(iii)] $\gcd (d_1, d_2)=1$,
\item[(iv)] for all odd primes $q \mid d_1$ we have $\; \displaystyle \kro{-d_2}{q}=1$,
\item[(v)] if $2 \mid d_1$ then $d_2 \equiv 7 \pmod{8}$.
\end{enumerate}
Moreover there are integers $s,~t$ such that 
\[
x=d_1 t, \qquad y=\mathrm{rad}(d_1) s,
\]
where $\mathrm{rad}(d_1)$ denotes the radical of $d_1$, and
\begin{equation}\label{eqn:st}
t^2+d_2=e s^p, \quad \gcd(t,d_2)=1, \quad s \neq 0 , 
% \qquad t \neq \pm 1, \quad \mathrm{gcd}(t,d_2)=1, \quad s \neq 0
\end{equation}
where
\begin{equation}\label{eqn:e}
e = \prod_{\substack{q\, \text{prime}\\q \mid d_1}} q^{p -2 v_q(d_1)} ,
\end{equation}
and $\; \mathrm{rad} (e)=\mathrm{rad}(d_1)$.
\end{lem}
\begin{proof}
	Suppose $(x,y,p)$ is a solution to equation~(\ref{eqn:p}) such
	that $y \neq 0$ and condition~(\ref{eqn:cond}) is satisfied.
	It is straightforward to see that condition~(\ref{eqn:cond}) 
	forces $\; \gcd(x^2,D)\;$ to be a square, say $d_1^2$ with
	$d_1 > 0$. We can therefore write $x=d_1 t$ and $D=d_1^2 d_2$
	for some integers $t,d_2$. Moreover, since 
	\[
		d_1^2=\gcd(x^2,D)=\gcd(d_1^2 t^2 , d_1^2 d_2)
		=d_1^2 \gcd(t^2 ,d_2),
	\]
	we see that $\gcd(t,d_2)=1$. Removing the common factors
	from $x^2+D=y^p$ we obtain $t^2+d_2=e s^p$ where $e$ is given
	by~(\ref{eqn:e}). The integrality of $e$ follows from the 
	condition~(\ref{eqn:cond}), and so does the equality of the 
	radicals $\mathrm{rad}(e)=\mathrm{rad}(d_1)$. 
	Note that (iii) follows from this equality of the radicals
	and the fact that $t,~d_2$ are coprime.
	We have thus proven (i), (ii), (iii) and it is now 
	easy to deduce (iv) and (v).
	Finally, the condition $\, s \neq 0$
	follows from the condition $\, y \neq 0$.
\end{proof}

\begin{defn}
Suppose $D$ is a non-zero integer and $(x,y,p)$ is a solution to
equation~(\ref{eqn:p}) with $y \neq 0$ and $p$ satisfying~(\ref{eqn:cond}).
Let $d_1,~d_2$ be as in the above Lemma and its proof (thus $d_1 >0$ and $\gcd(x,D)=d_1^2$
and $d_2=D/d_1^2$). We call the pair $(d_1,d_2)$ the signature of 
the solution $(x,y,p)$. We call the pair $(t,s)$ the simplification
of $(x,y)$ (or $(t,s,p)$ the simplification of $(x,y,p)$).
\end{defn}

In this terminology, Lemma~\ref{lem:remove} associates to any $D$ a finite
set of possible signatures $(d_1,d_2)$ for the solutions $(x,y,p)$ of 
equation~(\ref{eqn:p}) satisfying~(\ref{eqn:cond}) and $y \neq 0$.
To solve~(\ref{eqn:p}) it is sufficient to solve it under the assumption
that the solution's signature is $(d_1,d_2)$ for each possible signature.
\begin{ex}
For example, if $D=25$, there are two possible signatures satisfying
the conditions of Lemma~\ref{lem:remove}; these are
$(d_1,d_2)=(1,25)$ or $(5,1)$. If $(d_1,d_2)=(1,25)$, then
$x=t,~y=s$ and we must solve the equation
\[
t^2+25=s^p, \quad 5 \nmid t.
\]
However, if $(d_1,d_2)=(5,1)$, then $x=5t,~y=5 s$, and we must
solve the equation
\[
t^2+1=5^{p-2} s^p.
\]
In either case it is noted that the three terms of the resulting equation
are relatively coprime, which is important when we come to apply the modular
approach.
\end{ex}

\section{A Simplification of Cohn}\label{sec:Cohn}
We will soon apply our modular machinery to equations~(\ref{eqn:p})
with $D$ in the range~(\ref{eqn:Drange}). Before doing this
it is helpful to introduce a simplification due to Cohn
that will drastically reduce
the amount of computation needed later.
All the arguments presented in this
Section are found in Cohn's papers~\cite{Cohn1},~\cite{Cohn2}.
Cohn however assumes that $D \not \equiv 7 \pmod{8}$,
and the result that we state below is not formulated explicitly
that way in these papers. 
\begin{prop}\label{prop:Cohn}
Let $D=D_1^2 D_2$ where $D_2$ is square-free and $D_2 >0$. Suppose
that $(x,y,p)$ is a solution to equation~(\ref{eqn:p})
with $p$ satisfying~(\ref{eqn:cond}), and let $(d_1,d_2)$ 
be the signature of this solution. Then either
\begin{enumerate}
\item[(i)] $d_1 >1$,
\item[(ii)] or $\, D \equiv 7 \pmod 8$ and $2 \mid y$,
\item[(iii)] or $\, p$ divides the class number $h$ of the quadratic field 
$\Q(\sqrt{-D_2})$,
\item[(iv)] or $\, y=a^2 + D_2 b^2$ for some integers $a,~b$ such that
\[
b \mid D_1, \quad b \neq \pm D_1, \quad  p \mid (D_1^2 -b^2),
\]
and $a$ is a solution of the equation
\[
        \frac{1}{2 \sqrt{-D_2}} 
	\left[ (U+b \sqrt{-D_2})^p -(U-b \sqrt{-D_2})^p \right]=D_1,
\]
\item[(v)] or $D=1$, $(x,y)=(0,1)$, 
\item[(vi)] or $\, D_2 \equiv 3 \pmod{4}$ and $y=(a^2+D_2 b^2)/4$ for some odd integers $a,~b$ such that
\[
b \mid D_1, \quad  p \mid (4 D_1^2 -b^2),
\]
and $a$ is a solution of the equation
\[
        \frac{1}{2 \sqrt{-D_2}} \left[ (U+b \sqrt{-D_2})^p-(U-b \sqrt{-D_2})^p  \right]=2^p D_1.
\]
\end{enumerate}
\end{prop}
\begin{proof}
We only give a brief sketch. Suppose that (i), (ii), (iii) are false.
Then $(x+D_1 \sqrt{-D_2})=\alpha^p$ for some $\alpha$ in the ring of integers
of $\Q(\sqrt{-D_2})$. There are two possibilities. The first is
that $\alpha=a+b \sqrt{-D_2}$ for some integers $a$, $b$. 
By equating the imaginary parts
we deduce all of (iv) if $b \neq \pm D_1$. Thus suppose that  
$b= \pm D_1$. Letting $\beta=a - b\sqrt{-D_2}$
we see that
\[
\frac{\alpha^p -\beta^p}{\alpha-\beta}=\pm 1.
\]
If $\alpha/\beta$ is not a root of unity, then the 
left-hand side is the $p$-th term of a Lucas
sequence (with $p \geq 7$) and a deep Theorem of Bilu, Hanrot and Voutier
\cite{BHV} on primitive divisors of Lucas and Lehmer sequences
immediately gives a contradiction. Thus $\alpha/\beta$
is a root of unity and so equal to $\pm 1$, $\pm i$,
or $(\pm 1 \pm \sqrt{-3})/2$. Each case turns out to be
impossible, except for $\alpha=-\beta$ which together
with $b = \pm D_1$ implies (v). 

The second possibility for $\alpha$ is that $\alpha=(a+b\sqrt{-D_2})/2$
with $a,~b$ odd integers (and $-D_2 \equiv 1 \pmod{4}$).
Now (vi) follows quickly by equating the imaginary parts of
$(x+D_1 \sqrt{-D_2})=\alpha^p$.
\end{proof}

\subsection{Results II}
\begin{cor}\label{cor:Cohn}
Suppose $D$ belongs to our range~(\ref{eqn:Drange}) and $(x,y,p)$ is
a solution to equation~(\ref{eqn:p}) with $p$ satisfying
the condition~(\ref{eqn:cond}). If the solution $(x,y,p)$ is
missing from our Tables, then either
$D \equiv 7 \pmod{8}$ and $2 \mid y$, or $d_1 > 1$ 
where $(d_1,d_2)$ is the signature of the solution.
\end{cor}
\begin{proof}
We apply Proposition~\ref{prop:Cohn}. Using a short
{\tt MAGMA} program we listed all solutions arising from
possibilities (iv)--(vi) of that Proposition with $1 \leq D \leq 100$.
The only ones found in our range are $(x,y,p)=(0,1,p)$ for $D=1$
and $(x,y,p)=(\pm 8, 2,7)$ for $D=64$
and these are certainly in the Tables.

To prove the Corollary we merely have to take
care of possibility (iii) of the Proposition.
For $1 \leq D \leq 100$, and primes $p$
satisfying~(\ref{eqn:cond}), the only
case when $p$ could possibly divide the class number
of $\Q(\sqrt{-D_2})$ is $p=7$ and $D=71$ (in which case
$h=7$). 
We solved the equation~$x^2+71=y^7$ by reducing to Thue
equations as in Section~\ref{sec:Thue}. It took {\tt pari/gp}
about 30 minutes to solve these Thue equations, and
we obtained that the only solutions are $(x,y)=(\pm 46,3)$,
again in our Tables.
\end{proof}

\section{Level-Lowering}\label{sec:LL}
In this section we apply the modular approach to equation~(\ref{eqn:st})
under suitable, but mild, hypotheses.
Ordinarily, one would have to construct a Frey curve or curves
associated to our equation, show that the Galois representation is
irreducible (under suitable hypotheses) using the results of 
Mazur and others \cite{Mazur}, and modular by the work of Wiles and others 
\cite{Wi}, \cite{TW}, \cite{Mod}, and finally apply
Ribet's level-lowering Theorem \cite{Ribet}. Fortunately 
we are saved much trouble by the excellent paper of Bennett and Skinner \cite{BS},
which does all of this for equations of the form $Ax^n+By^n=Cz^2$;
it is noted that equation~(\ref{eqn:st}) is indeed of this form.

Let $D$ be a non-zero integer.
We shall apply the modular approach to the Diophantine equation
\begin{equation}\label{eqn:grn}
x^2+D=y^p, \qquad 
\quad x^2 \nmid D, \quad y \neq 0, \quad \text{and $p \geq 3$ is prime}
\end{equation}
or the equivalent equation for the simplification $(s,t)$
\begin{equation}\label{eqn:st1}
t^2+d_2=e s^p, 
 \qquad t \neq \pm 1, \quad \gcd(t,d_2)=1, \quad s \neq 0
\end{equation}
under the additional assumption that $p$ satisfies~(\ref{eqn:cond}).
%\begin{equation}\label{eqn:cond}
%\text{for all primes $q \mid D$ we have }
%p \geq v_q(D)+1.  \\
%\end{equation}
The assumptions made about $s$, $t$ %%!
in~(\ref{eqn:st1})
are there to ensure the non-singularity of the Frey curves, and the
absence of complex multiplication when we come to apply the modular
approach later on.
Before going on we note the following Lemma which in effect
says that there is no harm in making these additional assumptions for $D$ 
in our range~(\ref{eqn:Drange}).

\begin{lem}
There are no solutions to the equation~(\ref{eqn:p})
for $D$ in the range~(\ref{eqn:Drange})
with $y=0$, or $x^2 \mid D$,
except those listed in the Tables at the end.
\end{lem}
\begin{proof}
Clearly $y \neq 0$. 
We produced our list of solutions with $x^2 \mid D$ 
using a short {\tt MAGMA} program.
\end{proof}

Lemma~\ref{lem:remove} associates to each equation of the form~(\ref{eqn:grn})
finitely many signatures $(d_1,d_2)$ satisfying conditions~(i)--(v), %%!
and corresponding equations~(\ref{eqn:st}).
Following Bennett and Skinner \cite{BS} we associate a Frey curve $E_t$ 
to any potential solution of
equation~(\ref{eqn:st1}) according to 
Tables~\ref{table:LL1},~\ref{table:LL2},~\ref{table:LL3}.

\begin{table}
\caption{Frey Curves with $d_1$, $d_2$ odd.} %%!
\begin{tabular}{||c|c|c|c|c||}
\hline \hline
Case & Condition on $d_2$ & Condition on $t$ & Frey Curve $E_t$ & $L$ \\ 
\hline \hline
(a) & $d_2 \equiv 1 \pmod{4}$  & & $Y^2=X^3+2 tX^2 -d_2 X$ & $2^5 $ \\
\hline
(b) & $d_2 \equiv 3 \pmod{8}$ & & $Y^2=X^3+2 tX^2 +(t^2+d_2) X$ & $2^5 $ \\
\hline
(c) & $d_2 \equiv 7 \pmod{8}$ & $t$ even & $Y^2=X^3+2 tX^2 +(t^2+d_2) X$ & $2^5 $ \\
\hline
(d) & $d_2 \equiv 7 \pmod{8}$ & $t \equiv 1 \pmod{4}$ & 
$Y^2+XY=X^3 + \left( \frac{t-1}{4}\right) X^2 + \left( \frac{t^2+d_2}{64}\right)X$ & $2 $\\
\hline\hline
\end{tabular}
\label{table:LL1}
\end{table}

\medskip

\begin{table}
\caption{Frey Curves with $d_1$ even, $d_2$ odd.}
\begin{tabular}{||c|c|c|c||}
\hline \hline
Case & Conditions on $t,s,p$ & Frey Curve $E_t$ & $L$ \\
\hline \hline
(e) & $t \equiv 1 \pmod{4}$  &  
	$Y^2+XY=X^3+\left(\frac{t-1}{4} \right)X^2+\left(\frac{t^2+d_2}{64} \right)X$ & $1$\\
%		& either $p \geq v_2(D)+7$ or $s$ even	&	& \\
%\hline
%(f) & $t \equiv 1 \pmod{4}$,  &
%        $Y^2=X^3+ t X^2+\left(\frac{t^2+d_2}{4}\right)X$ & $2^2$\\
%	        & $p = v_2(D)+5$ and $s$ odd &       & \\
%\hline
%(g) & $t \equiv 1 \pmod{4}$, &
%        $Y^2=X^3+t X^2+\left(\frac{t^2+d_2}{4} \right)X$ & $2^4$\\
%	        & $p = v_2(D)+3$ and $s$ odd &       & \\
%\hline
%(h) &  $t \equiv 1 \pmod{4},$&
%        $Y^2+XY=X^3+ 2 t X^2+(t^2+d_2) X$ & $2^6$\\
%	         & $p =v_2(D)+1$ and $s$ odd &       & \\
\hline\hline
\end{tabular}
\label{table:LL2}
\end{table}

\medskip

\begin{table}
\caption{Frey Curves with $d_1$ odd, $d_2$ even.}
\begin{tabular}{||c|c|c|c|c||}
\hline \hline
Case & Condition on $d_2$ & Condition on $t$ & Frey Curve $E_t$ & $L$ \\
\hline \hline
(f) & $v_2(d_2)=1$ & & $Y^2=X^3+2t X^2 -d_2 X$ & $ 2^6 $ \\
\hline
(g) & $d_2 \equiv 4 \pmod{16}$ & $t \equiv 1 \pmod{4}$ & 
	$Y^2=X^3+tX^2-\frac{d_2}{4}X$ & $2 $\\
\hline
(h) & $d_2 \equiv 12 \pmod{16}$ & $t \equiv 3 \pmod{4}$ &
       $Y^2=X^3+tX^2-\frac{d_2}{4}X$ & $2^2 $\\
\hline
(i) & $v_2(d_2)=3$ & $t \equiv 1 \pmod{4}$ & $Y^2=X^3+t X^2 -\frac{d_2}{4}X$ & $2^4 $\\
\hline
(j) & $v_2(d_2)=4,5$ & $t \equiv 1 \pmod{4}$ & $Y^2=X^3+t X^2 -\frac{d_2}{4}X$ & $2^2 $\\
\hline
(k) & $v_2(d_2) = 6$ & $t \equiv 1 \pmod{4}$ & 
$Y^2+XY=X^3+ \left( \frac{t-1}{4} \right)X^2 - \frac{d_2}{64} X$ & $2^{-1} $\\
\hline
(l) & $v_2(d_2) \geq 7$ & $t \equiv 1 \pmod{4}$ &
$Y^2+XY=X^3+ \left( \frac{t-1}{4} \right)X^2 - \frac{d_2}{64} X$ & $1$\\
\hline \hline
\end{tabular}
\label{table:LL3}
\end{table}

The three tables are divided into cases (a)--(l).
We know that $d_1,~d_2$ are coprime, and hence at most one of them is even.
The possibility that $d_1,~d_2$ are both odd is dealt with in Table~\ref{table:LL1}. 
In cases (a),~(b), a simple modulo $8$ argument convinces us that $t$ is odd.
However for cases (c) and (d) -- where $d_1$ is odd and $d_2 \equiv 7 \pmod{8}$ --
the integer $t$ can be either odd or even and we assign different Frey curves for each possibility.
When $t$ is odd (case (d)) we add the assumption that $t \equiv 1 \pmod{4}$. 
This can be achieved by interchanging $t$ with $-t$ if necessary.

Table~\ref{table:LL2} deals with the possibility of even $d_1$, and Table~\ref{table:LL3} 
with the possibility of even $d_2$. In both these cases $t$ is necessarily odd, 
and the congruence condition on
$t$ can again be achieved by interchanging $t$ with $-t$ if necessary. 

\begin{prop}\label{prop:LL}
Suppose $D$, $d_1$, $d_2$ %%!!
are non-zero integers that satisfy~(i)--(v) of Lemma~\ref{lem:remove}.
Suppose also
$p$ is a prime number satisfying the condition~(\ref{eqn:cond}),
and let $e$ be as defined in~(\ref{eqn:e}). Suppose
that $(t,s)$ is a solution of
equation~(\ref{eqn:st1}) and satisfying the supplementary
condition (if any) on $t$ in Tables~\ref{table:LL1},~\ref{table:LL2},~\ref{table:LL3}. 
Let $E_t$ and $L$ be as in these tables, 
and write $\rho_p(E_t)$ be the Galois
representation on the $p$-torsion of $E_t$.
Then the representation $\rho_p(E_t)$ arises from a cuspidal newform of weight $2$
and level $N=L \, \mathrm{rad}(D)$. 
\end{prop}
\begin{proof} 
The paper of Bennett and Skinner \cite{BS}
gives an exhaustive recipe for Frey curves and level-lowering for 
equations of the form $A x^n + B y^n= C z^2$ under the assumption
that the three terms in the equation are coprime. After a little 
relabeling, their results apply to our equation~(\ref{eqn:st1})
and the Lemma follows from Sections~2,~3 of their paper. It is here
that we need the assumptions $t \neq \pm 1$ and $s \neq 0$
made in~(\ref{eqn:st1})
\end{proof}

It is convenient to indulge in the following abuse of language.
\begin{defn}
If $(t,s,p)$ is a solution to equation~(\ref{eqn:st1}) 
and if the representation $\rho_p(E_t)$ arises from
a cuspidal newform $f$, then we say that solution $(t,s,p)$ arises
from the newform $f$ (via the Frey curve $E_t$), or that the newform $f$ 
gives rise to the solution $(t,s,p)$. If $(t,s,p)$ is the 
simplification of $(x,y,p)$ then we say that $(x,y,p)$ arises
from the newform $f$.

If the newform $f$ is rational, and so corresponds to an elliptic curve $E$,
then we also say that the solution $(t,s,p)$ (or $(x,y,p)$) arises from $E$.
\end{defn}

\subsection{A Summary}
It may be helpful for the reader to summarize what we have done and 
where we are going. Given a non-zero integer $D$ we would like 
to solve equation~(\ref{eqn:grn}). We can certainly write down all
solutions with $y=0$ or with $x^2 \mid D$. We can also solve (at least in
principle) all cases where $p$ violates 
condition~(\ref{eqn:cond}) by reducing to Thue equations as in Section~\ref{sec:Thue}.
We can thus reduce to 
equation~(\ref{eqn:grn}) and assume that $p$ satisfies condition~(\ref{eqn:cond}).

Next, we can write down a list of signatures $(d_1,d_2)$ 
satisfying conditions~(i)--(v) of Lemma~\ref{lem:remove}.
We reduce the solution of equation~(\ref{eqn:grn}) to solving
for each signature $(d_1,d_2)$ the equation~(\ref{eqn:st1}).
Now we associate to the signature $(d_1,d_2)$ one or more Frey curves $E_t$
and levels $L$, so that any solution to~(\ref{eqn:st1}) arises from
some newform $f$ at level $L$ via the Frey curve $E_t$.

Finally (and this is to come) we must show how to solve~(\ref{eqn:st1})
under the assumption that the solution arises from a newform $f$ via
a Frey curve $E_t$. If we can do this for each newform $f$ and Frey
curve $E_t$ then we will have completed the solution of our 
equation~(\ref{eqn:p}).

As we shall see, the assumption that a solution arises from a particular
newform is a very strong one, for it imposes congruence conditions on
$\, t \,$ modulo all but finitely many primes $\, l$.

\subsection{Congruences}
For an elliptic curve $E$ we write $\sharp E(\F_l)$ for the number of points
on $E$ over the finite field $\F_l$, and let $a_l(E)=l+1- \sharp E(\F_l)$.
\begin{lem}\label{lem:cong1}
With notation as above, suppose that the Galois representation $\rho_p(E_t)$
arises from a cuspidal newform with Fourier expansion around infinity 
\begin{equation}\label{eqn:Fourier}
f=q+\sum_{n \geq 2} c_n q^n,
\end{equation}
of level $N$ (given by Proposition~\ref{prop:LL}) and defined over a number field $K/\Q$.
Then there is a place $\Pp$ of $K$ above $p$ such that for every prime
$l \nmid 2pD$ we have
\[
\begin{array}{ll}
a_l(E_t) \equiv c_l \pmod{\Pp} & \text{if $t^2+d_2 \not \equiv 0 \pmod{l}$ (or equivalently $l \nmid s$)}, \\
l+1 \equiv \pm c_l \pmod{\Pp} & \text{if $t^2+d_2  \equiv 0 \pmod{l}$ (or equivalently $l \mid s$)}.  \\
\end{array}
\]
\end{lem}
\begin{proof}
The Lemma is standard (see \cite[page 196]{Se87}, \cite[page 7]{BS},
\cite[Proposition 5.4]{Kra}, etc.). The conditions $l \nmid 2D$ and 
$l \nmid s$ together imply that $l$ is a prime of good reduction for $E_t$,
whereas the conditions $l \nmid 2D$ and $l \mid s$ imply
that $l$ is a prime of multiplicative reduction. 
\end{proof}

When the newform $f$ is rational, there is an elliptic curve $E$  
defined over $\Q$ whose conductor is equal to the level of the newform $f$
such that $a_l(E)=c_l$ for all primes $l$. 
In this case we can be a little more precise than in 
Lemma~\ref{lem:cong1}, thanks to a result of Kraus and Oesterl\'{e}.

\begin{lem}\label{lem:cong2}
With notation as above, suppose that the Galois representation $\rho_p(E_t)$
arises from a rational cuspidal newform $f$ corresponding to an elliptic curve $E/\Q$.
Then for all primes $l \nmid 2D$ we have
\[
\begin{array}{ll}
a_l(E_t) \equiv a_l(E) \pmod{p} & \text{if $t^2+d_2 \not \equiv 0 \pmod{l}$ 
(or equivalently $l \nmid s$)},
\\
l+1 \equiv \pm a_l(E) \pmod{p} & \text{if $t^2+d_2  \equiv 0 \pmod{l}$ 
(or equivalently $l \mid s$)}.
\\
\end{array}
\]
\end{lem}
\begin{proof}
This Lemma does appear to be a special case of Lemma~\ref{lem:cong1}; however we
do allow in this Lemma the case $l=p$ which was excluded before. In fact 
Lemma~\ref{lem:cong1} together with a result of 
Kraus and Oesterl\'{e} \cite[Proposition 3]{KO} implies that the representations $\rho_p(E_t)$
and $\rho_p(E)$ are semi-simply isomorphic. In this case the result of Kraus and Oesterl\'{e}
also tells us that $a_l(E_t) \equiv a_l(E) \pmod{p}$ if the prime $l$ is a prime of good reduction
for both curves, and $a_l(E_t)a_l(E) \equiv l+1 \pmod{p}$ if $l$ is a prime of good reduction for
one of them and a prime of multiplicative reduction for the other. Now since $ l \nmid 2D$
we see that $l$ does not divide the conductor $N$ of $E$ (which is also the
level of the newform $f$ as given by Proposition~\ref{prop:LL}). If $l \mid s$ %%!
then $l$ is a prime of multiplicative reduction for $E_t$ and then $a_l(E_t)=\pm 1$.
The Lemma follows.   
\end{proof}

%\section{Absence of Newforms}
%Lemma~\ref{lem:remove} and Proposition~\ref{prop:LL}
%lead us to associate solutions to equation~(\ref{eqn:grn})
%with $p$ satisfying condition~(\ref{eqn:cond}),
%with newforms of certain levels. If there are 
%no newforms of the predicted levels, we immediately deduce that there
%are no solutions to equation~(\ref{eqn:grn}).
%With the help of a {\tt MAGMA} program we find all values of $D$ in
%the range $-100 \leq D \leq 100$ where there are no newforms  
%at the predicted levels, and deduce the following result.
%\begin{cor}
%Let $D$ be an integer belonging to the set
%\[
%\left\{ 
%-64,-60,-32,-28,-16,-12,-4,4,16,32,36,64
%\right\}
%\]
%Then the equation~(\ref{eqn:grn}) does not have any solutions with 
%$p$ satisfying condition~(\ref{eqn:cond}).
%\end{cor}

\section{Eliminating Exponents: Method I}

We now focus on equations of the form~(\ref{eqn:st1}) where, as always, $p$ satisfies
~(\ref{eqn:cond}). Proposition~\ref{prop:LL} tells us that if $(t,s,p)$ is a 
solution to~(\ref{eqn:st1}), then it arises from a newform
of a certain level (or levels) and all these can be determined. Let us say
that these newforms are $f_1,\ldots,f_n$. Then to solve equation~(\ref{eqn:st1})
it is sufficient to solve it, for each $i$, under the assumption that the solution arises from
the newform $f_i$. 
We give three methods for attacking equation~(\ref{eqn:st1}) under the assumption that
the solution arises from a particular newform $f$. 
%\subsection{Method I}

If successful, the first method will prove that the equation~(\ref{eqn:st1})
has no solutions except possibly for finitely many exponents $p$ and these are determined 
by the method. This method is actually quite standard. We believe that the basic idea 
is originally due to Serre \cite[pages 203--204]{Se87}.
It is also found in Bennett and Skinner  \cite[Proposition 4.3]{BS}.
We shall however give a more careful version
than is found in the literature, thereby maximizing the probability of success.

\begin{prop}\label{prop:methodI}
{\rm (Method I)}
 Let $D$, $d_1$, $d_2$ be a %%!
 triple of integers satisfying (i)--(v) of Lemma~\ref{lem:remove}.
% Let $E_t$ be a Frey curve as in Table~\ref{table:LL}, and let $N=L \mathrm{rad}(D)$
%where $L$ is given as in the last column of the Table. 
Let $f$ be a newform 
with Fourier expansion as in~(\ref{eqn:Fourier}) having coefficients in the 
ring of integers of a number field $K$, and let $\N_{K/\Q}$ denote the norm map. 
% Suppose that the newform $f$ gives rise to the Galois representation $\rho_p(E_t)$. 
If $l\nmid 2D$ is prime, let
\[
        B_l^{\prime \prime}(f) =
        \mathrm{lcm} \left\{ \N_{K/\Q} (a_l(E_t) -c_l) \, \,
                                : \quad t \in \F_l, \quad t^2+d_2 \not \equiv 0 \pmod{l} \right\},
\]
\begin{equation*}
        B_l^{\prime}(f) = 
	\begin{cases}
                        B_l^{\prime \prime}(f) & \text{if\ } \kro{-d_2}{l}=-1, \\
                        \mathrm{lcm} \left\{ B_l^{\prime \prime}(f) ,  \,
                                        \N_{K/\Q} (l+1+c_l), \,
                                        \N_{K/\Q} (l+1 -c_l) \right\} &
                                        \text{if\ } \kro{-d_2}{l}=1,
	\end{cases}
\end{equation*}
and
\begin{equation*}
	B_l(f) = 
	\begin{cases}
        l\, B_l^\prime (f) & \text{if $K \neq \Q$}, \\
        B_l^\prime (f) & \text{if $K = \Q$}.
        \end{cases}
\end{equation*}
If $p$ satisfies condition~(\ref{eqn:cond}), and if 
$(t,s,p)$ is a solution to equation~(\ref{eqn:st1}) 
arising from the newform $f$
% with $t$ satisfying the supplementary  condition for $E_t$ in Table~\ref{table:LL}, 
then $p \mid B_l(f)$.
\end{prop}
\begin{proof}
The Proposition follows almost immediately from 
Lemmas~\ref{lem:cong1},~\ref{lem:cong2}.
\end{proof}

Under the assumptions made (in this Proposition), Method I eliminates all but finitely
many exponents $p$ provided
of course that the integer $B_l(f)$ is non-zero. 
Accordingly, we shall say that Method I is successful if there exists some prime $l \nmid 2D$
so that $B_l(f) \neq 0$. There are two situations where Method I is guaranteed to succeed:
\begin{itemize}
\item If the newform $f$ is not rational. In this case, for infinitely many primes $l$,
the Fourier coefficient $c_l \not \in \Q$ and so all the 
differences $a_l(E_t)-c_l$ and $l+1 - c_l$ are 
certainly non-zero, immediately implying that $B_l(f) \neq 0$.

\item Suppose that the newform $f$ is rational, and so corresponds to an elliptic
curve $E$ defined over $\Q$. Suppose that $E$ has no non-trivial $2$-torsion.
By the Chebotarev Density Theorem we know that $l+1 - a_l(E)= \sharp E(\F_l)$
is odd for infinitely many primes $l$. Let $l \nmid 2D$ be any such prime.
>From the models for the Frey curves in Tables~\ref{table:LL1},~\ref{table:LL2},~\ref{table:LL3}
we see that the Frey curve $E_t$ has non-trivial $2$-torsion, and so $l+1 -a_l(E_t)
= \sharp E_t(\F_l)$ is even for any value of $t \in \F_l,\, t^2+d_2 \neq 0$. In this case 
$a_l(E_t)-c_l=a_l(E_t)-a_l(E)$ must be odd and cannot be zero. Similarly, the Hasse--Weil
bound $\vert c_l \vert \leq 2 \sqrt{l}$ implies that $l+1 \pm c_l \neq 0$. Thus $B_l(f)$
is non-zero in this case and Method I is successful.
\end{itemize}

\section{Eliminating Exponents: Method II}
The second method is adapted from the ideas of Kraus \cite{Kra} (see also \cite{SC}).
It can only be applied to one prime (exponent) $p$ at a time, and if successful
it does show that there are no solutions to~(\ref{eqn:st1}) for that particular exponent.

Let us briefly explain the idea of this second method. Suppose 
$f$ is a newform with Fourier expansion as in~(\ref{eqn:Fourier}),
and suppose $p \geq 7$ is a prime. We are interested in solutions to equation~(\ref{eqn:st1})
arising from $f$. Choose a small integer
$n$ so that $l=np+1$ is prime with $l \nmid D$. Suppose $(t,s)$ is a solution to
equation~(\ref{eqn:st1}) arising from newform $f$. 
Then working modulo $l$ we see that $d_1^2 t^2+D=y^p$  %$%
is either $0$ or an $n$-th root of unity. Since $n$ is small we can list
all such $t$ in $\F_l$, and compute $c_l$ and $a_l(E_t)$ for each $t$ in our 
list. We may then find that for no $t$ in our list are the relations in 
Lemma~\ref{lem:cong1} satisfied. In this case we have a contradiction, and we deduce
that the are no solutions to equation~(\ref{eqn:st1}) arising from $f$ for our particular exponent $p$.

Let us now write this formally. Suppose $p \geq 7$ is a prime number, and $n$
an integer such that $l=np+1$ is also prime and $l \nmid D$. Let
\[
\mu_n(\F_l)=\left\{ \zeta \in \F_l^* \quad : \quad \zeta^n=1 \right\}.
\]
Define 
\[
A(n,l)= \left\{
	\zeta \in \mu_n(\F_l) \quad : \quad \left( \frac{\zeta-D}{l}\right)= 0 \text{\ or\ } 1
	\right\}.
\]
For each $\zeta \in A(n,l)$, let $\delta_\zeta$ be an integer satisfying
\[
\delta_\zeta^2 \equiv (\zeta-D)/d_1^2 \pmod{l}.
\]
It is convenient to write $a_l(\zeta)$ for $a_l(E_{\delta_\zeta})$.

We now give our sufficient condition for the insolubility of~(\ref{eqn:st1}) for the
given exponent $p$.

\begin{prop}\label{prop:methodII}
{\rm(Method II)} Let $D$, $d_1$, $d_2$ be a triple of integers satisfying (i)--(v) of Lemma~\ref{lem:remove},
and $p \geq 7$ be a prime satisfying condition~(\ref{eqn:cond}).
% Let $E_t$ be a Frey curve as in Table~\ref{table:LL}, and let $N=L \, \mathrm{rad}(D)$
% where $L$ is given as in the last column of the Table. 
Let $f$ be a newform 
with Fourier expansion as in~(\ref{eqn:Fourier}) defined over a number field $K$.
Suppose there exists an integer $n \geq 2$ satisfying
the following conditions:
\begin{enumerate}
\item[(a)] The integer $l=np+1$ is prime, and $l \nmid D$.
\item[(b)] Either $\displaystyle\left(\frac{-d_2}{l}\right) = -1$, or 
$\displaystyle \, \, p \nmid \N_{K/\Q} (4-c_l^2)$.
\item[(c)] For all $\zeta \in A(n,l)$ we have 
\[
\left\{
\begin{array}{ll}
p \nmid \N_{K/\Q} (a_l(\zeta)-c_l) & \text{if $l \equiv 1 \pmod{4}$,} \\
p \nmid \N_{K/\Q} (a_l(\zeta)^2-c_l^2) & \text{if $l \equiv 3 \pmod{4}$.}
\end{array}
\right.
\]
\end{enumerate}
Then the equation~(\ref{eqn:st1}) does not have any solutions for the given exponent $p$
arising from the newform $f$. 
\end{prop}
\begin{proof}
Suppose that the hypotheses of the Proposition are satisfied, and that
$(t,s)$ is a solution to equation~(\ref{eqn:st1}).

First we show that $t^2+d_2 \not \equiv 0 \pmod{l}$. Suppose otherwise. Thus
$t^2+d_2 \equiv 0 \pmod{l}$ and so $l \mid s$.
In this case $\kro{-d_2}{l}=1$, and from (b) we know that 
$\, p \nmid \N_{K/\Q} (4-c_l^2)$.
However, by Lemma~\ref{lem:cong1} we know that $\pm c_l \equiv l+1 \equiv 2 \pmod{\Pp}$
for some place $\Pp$ of $K$ above $p$, and we obtain a contradiction showing that 
$t^2+d_2 \not \equiv 0 \pmod{l}$. 

>From equation~(\ref{eqn:st1}) and the definition of $e$ in~(\ref{eqn:e}) 
we see the existence of some $\zeta \in A(n,l)$ such that 
\[
d_1^2 t^2+D \equiv \zeta \pmod{l} \quad \text{and} \quad t \equiv \pm \delta_\zeta \pmod{l}.
\] 
Replacing $t$ by $-t$ in the Frey curve $E_t$ has the effect of twisting the curve
by $-1$ (this can be easily verified for each Frey curve in 
Tables~\ref{table:LL1},~\ref{table:LL2},~\ref{table:LL3}).
Thus $a_l(\zeta) = a_l(E_t)$ if $l \equiv 1 \pmod{4}$ and 
$a_l(\zeta)=\pm a_l(E_t)$ if $l \equiv 3 \pmod{4}$. Moreover, by Lemma~\ref{lem:cong1} 
we know that $a_l(E_t) \equiv c_l \pmod{\Pp}$ for some place $\Pp$ of $K$
above $p$. This clearly contradicts (c).
Hence there is no solution to~(\ref{eqn:st1}) arising from $f$ for the given exponent $p$.
\end{proof}

If the newform $f$ is rational and moreover corresponds to
an elliptic curve with $2$-torsion, then it is possible to
strengthen the conclusion of Proposition~\ref{prop:methodII}
by slightly strengthening the hypotheses. The following
variant is far less costly in computational terms as we explain below.

\begin{prop}\label{prop:methodIIrational}
{\rm(Method II)} Let $D$, $d_1$, $d_2$ be a triple of integers satisfying (i)--(v) of Lemma~\ref{lem:remove},
and $p$ be a prime satisfying condition~(\ref{eqn:cond}).
Let $f$ be a rational newform corresponding to elliptic curve $E/\Q$ with
$2$-torsion.
Suppose there exists an integer $n \geq 2$ satisfying
the following conditions:
\begin{enumerate}
\item[(a)] The integer $l=np+1$ is prime, $\displaystyle l \leq \frac{p^2}{4}$ and $l \nmid D$.
\item[(b)] Either $\displaystyle \left(\frac{-d_2}{l}\right) = -1$, or
$\displaystyle \; a_l(E)^2 \not \equiv 4 \pmod{p}$.
\item[(c)] For all $\zeta \in A(n,l)$ we have
\[
\left\{
\begin{array}{ll}
a_l(\zeta) \neq a_l(E) & \text{if $l \equiv 1 \pmod{4}$,} 
\\
a_l(\zeta) \neq \pm a_l(E) & \text{if $l \equiv 3 \pmod{4}$.}
\end{array}
\right.
\]
\end{enumerate}
Then the equation~(\ref{eqn:st1}) does not have any solutions for the given exponent $p$
arising from the newform $f$.
\end{prop}
\begin{proof}
Comparing this with Proposition~\ref{prop:methodII} we see that
it is sufficient to show, under the additional assumptions,
that if $a_l(\zeta)^2 \equiv a_l(E)^2 \pmod{p}$ then
$a_l(\zeta)= \pm a_l(E)$, and if $a_l(\zeta) \equiv a_l(E) \pmod{p}$ then
$a_l(\zeta)= a_l(E)$.

Suppose that $a_l(\zeta)^2 \equiv a_l(E)^2 \pmod{p}$ (the other case
is similar).  Hence $a_l(\zeta)\equiv \pm a_l(E) \pmod{p}$. Now note that both
elliptic curves under consideration here have $2$-torsion. Hence
we can write
\[
a_l(\zeta) = 2 b_1 \quad \text{and} \quad a_l(E)= 2 b_2
\]
for some integers $b_1,~b_2$. Moreover, by the Hasse--Weil
bound we know that $|b_i| \leq  \sqrt{l}$.
Thus
\[
b_1 \equiv \pm b_2 \pmod{p} \quad \text{and} \quad |b_1+b_2|,~|b_1-b_2| \leq 2 \sqrt{l} < p
\]
since $\displaystyle l< \frac{p^2}{4}$. Thus $b_1= \pm b_2$ and this completes the proof.
\end{proof}

It remains to explain how this improves our computation. To apply
Proposition~\ref{prop:methodII} for a particular prime $p$ we need
to find a prime $l$ satisfying conditions (a), (b), (c). The
computationally expensive part is to compute $a_l(E)=c_l$
and $a_l(\zeta)$ for all $\zeta \in A(n,l)$.
Let us however consider the application of Proposition~\ref{prop:methodIIrational}
rather than Proposition~\ref{prop:methodII}.
The computation proceeds as before by checking conditions (a), (b) first.
When we come to condition (c), we note that what we have to check
is that
\[
\left\{
\begin{array}{ll}
E_\zeta(\F_l) \neq l+1-a_l(E) & \text{if $l \equiv 1 \pmod{4}$,} \\
E_\zeta(\F_l) \neq l+1 \pm a_l(E) & \text{if $l \equiv 3 \pmod{4}$,} %%!
\end{array}
\right.
\]
for each $\zeta \in A(n,q)$. Rather than compute $a_l(\zeta)$
for each such $\zeta$, we first pick a random point in $E_\zeta(\F_l)$,
and check whether it is annihilated by $l+1-a_l(E)$ if
$p \equiv 1 \pmod{4}$ and either of the integers $l+1 \pm a_l(E)$ if $p \equiv 3 \pmod{4}$.
Only if this is the case do we need to compute $a_l(\zeta)$ to test
condition (c) in the Proposition. In practice, for primes $p$ of about $10^9$, this
brings a $10$-fold speed-up in program run time for Method II.

\section{Eliminating Exponents: Method III}\label{sec:methodIII}

 Occasionally, Methods I and II fail to establish the non-existence of solutions
 to an equation of the form~(\ref{eqn:st1}) for a particular exponent $p$
 even when it does seem that this equation has no solutions. 
 The reasons for this failure are not clear to us. We shall however give a third
 method, rather similar in spirit to Kraus' method (Method II), but requiring
 stronger global information furnished by Proposition~\ref{prop:Thue1}.
%It is appropriate to start with an informal explanation of the method.
%We know from Proposition~\ref{prop:Thue1} that if $(x,y)$
%is a solution to equation~(\ref{eqn:rn}) then
%$x+D_1 \sqrt{-D_2} = \gamma \beta^p$, for some
%$\gamma \in \OO$ and $\beta \in \OO$. Now applying the modular
%ideas, $x$ can be predicted to belong to a small set modulo $l$.
%However, if $l=np+1$ is a prime which splits in 
%$\Q(\sqrt{-D_2})$, and $n$ is small 
%then $\beta^p$ belongs to a small subset modulo $l$.

Suppose that $D$, $d_1$, $d_2$ are integers satisfying conditions
(i)--(v) of Lemma~\ref{lem:remove}. Let $E_t$ be one of the 
Frey curves associated to equation~(\ref{eqn:st1}), and 
let $f$ a newform of the level predicted by Proposition~\ref{prop:LL}
with Fourier expansion as in~(\ref{eqn:Fourier}),
defined over a number field $K$.
Define $\Tc_l(f)$ to be the set of $\tau \in \F_l$
such that
\begin{itemize}
\item either $p \mid \N_{K/\Q} (a_l(E_\tau)-c_l)$ and $\tau^2+d_2 \not \equiv 0 \pmod{l}$,
\item or $p \mid \N_{K/Q}(l+1 \pm c_l)$ and $\tau^2+d_2 \equiv 0 \pmod{l}$.
\end{itemize}
%As in Section~\ref{sec:Thue} we must consider separately the case where
%$-D$ is a non-square, and the case where $-D$ is a square.
We suppose that $-D$ is not a square and follow the notation of 
Section~\ref{sec:Thue}. Fix a prime $p$ satisfying~(\ref{eqn:cond}).
Suppose $l$ is a prime satisfying the following conditions:
\begin{enumerate}
\item[(a)] $l \nmid 2D$.
\item[(b)] $l=np+1$ for some integer $n$.
\item[(c)] $\displaystyle \kro{-D_2}{l}=1$. Thus $l$ splits in $\cL=\Q(\sqrt{-D_2})$ 
and we let $\gl_1$ and
$\gl_2$ be the prime ideals above $l$.
\item[(d)] Each $\gamma \in \Gamma$ is integral at $l$; what we mean by this
is that each $\gamma$ belongs to the intersection of the localizations 
$\OO_{\gl_1} \cap \OO_{\gl_2}$.
\end{enumerate}
We denote the two natural reduction maps by 
$\theta_1,\, \theta_2: \OO_{\gl_1} \cap \OO_{\gl_2} \rightarrow \F_l$.
These of course correspond to the two square-roots for $-D_2$
in $\F_l$, and are easy to compute. 

%We know from Proposition~\ref{prop:Thue1} that if $(x,y)$
%is a solution to equation~(\ref{eqn:rn}) then 
%$x+D_1 \sqrt{-D_2} = \gamma \beta^p$, for some
%$\gamma \in \OO$ and $\beta \in \OO$. 

%Now applying the modular
%ideas, $x$ can be predicted to belong to a small set modulo $l$.

Now let $\Gamma_l$ be the set of $\gamma \in \Gamma$ satisfying the
condition: there exists $\tau \in \Tc_l(f)$ such that
\begin{itemize}
\item $(d_1 \tau + D_1 \theta_1 ( \sqrt{-D_2}))^n \equiv \theta_1(\gamma)^n$ or $0 \pmod{l}$, and
\item  $(d_1 \tau + D_1 \theta_2(\sqrt{-D_2}))^n \equiv \theta_2(\gamma)^n$ or $0 \pmod{l}$.
\end{itemize}

\begin{prop}\label{prop:methodIII1}
{\rm (Method III)}
Let $p$  be a prime satisfying condition~(\ref{eqn:cond}).
Let $S$ be a set of primes $l$ satisfying the conditions (a)--(d) above.
With notation as above, if the newform $f$ gives rise to a solution $(t,s)$
to equation~(\ref{eqn:st1}), then $d_1 t + D_1 \sqrt{-D_2}=\gamma \beta^p$
for some $\beta \in \OO$ and some $\gamma \in \cap_{l \in S} \Gamma_l$.

In particular, if $\cap_{l \in S} \Gamma_l$ is empty, then the newform $f$
does not give rise to any solution to equation~(\ref{eqn:st1}) for
the given exponent $p$.
\end{prop}
\begin{proof}
Suppose that $(t,s)$ is a solution to equation~(\ref{eqn:st1}) arising from
newform $f$ via the Frey curve $E_t$. Clearly $\theta_1(t)=\theta_2(t)$
is simply the reduction of $t$ modulo $l$.
Let $\tau=\theta_1(t)=\theta_2(t) \in \F_l$.
It follows from Lemma~\ref{lem:cong1}
that $\tau \in \Tc_l(f)$. Let $(x,y)$ be the solution
to equation~(\ref{eqn:grn}) corresponding to $(t,s)$. Thus $x=d_1 t$.
We know by Proposition~\ref{prop:Thue1} that
\[
d_1 t+ D_1 \sqrt{-D_2} = \gamma \beta^p,
\]
for some $\gamma \in \Gamma$ and $\beta \in \OO$.
Applying $\theta_i$ to both sides and taking $n$-th powers (where
we recall that $l=np+1$) we obtain
\[
(d_1 \tau +D_1 \theta_i(\sqrt{-D_2}))^n \equiv \theta_i(\gamma)^n \theta_i(\beta)^{l-1} \pmod{l}.
\]
However $\theta_i(\beta)^{l-1} \equiv 0$ or $1 \pmod{l}$. Thus $\gamma \in \Gamma_l$
as defined above. The Proposition follows.
\end{proof}

%\subsection{The Case where $-D$ is a Square}
%We suppose that $D= -D_1^2$, and we follow the notation
%of Section~\ref{sec:Thue}. 
%Fix a prime $p \geq 7$. Suppose that
%$l$ is a prime satisfying the conditions.
%\begin{itemize}
%\item $l \nmid 2D$.
%\item $l=np+1$ for some integer $n$.
%\end{itemize}
%Define $\Gamma_l$ to be the set of $(\alpha,\beta) \in \Gamma$
%satisfying the condition: there exists $t \in \Tc_l(f)$ such that
%\begin{itemize}
%\item $(d_1 t + D_1)^n \equiv \alpha^n$ or $0 \pmod{l}$,
%\item and $(d_1 t - D_1)^n \equiv \beta^n$ or $0 \pmod{l}$.
%\end{itemize}
%
%\begin{prop}\label{prop:methodIII2}
%{\rm (Method III)}
%Let $p \geq 7$  be a prime satisfying condition~(\ref{eqn:cond}).
%Let $S$ be a set of primes $l$ satisfying the above conditions.
%With notation as above, if the newform $f$ gives rise to a solution $(t,s)$
%to equation~(\ref{eqn:st}), then 
%\[
%d_1 t= \frac{1}{2} \left(\alpha U^p + \beta V^p \right),
%\]
%for some pair $(\alpha,\beta) \in \cap_{l \in S} \Gamma_l$, and some
%integral solution $(U,V)$ to the Thue equation
%\[
%\alpha U^p - \beta V^p=2 D_1.
%\]
%
%In particular, if $\cap_{l \in S} \Gamma_l$ is empty, then the newform $f$
%does not give rise to any solution to equation~(\ref{eqn:st}) for
%the given exponent $p$.
%\end{prop}

\section{Examples}

It is clear that our three modular methods require
computations of newforms of a given level. Fortunately
the computer algebra suit {\tt MAGMA} has a package
completely devoted to such computations; the theory 
for these computations is explained by Cremona
\cite{Cre} for rational newforms, and by Stein \cite{Stein}
in the general case. 

\begin{ex} {\bf Absence of Newforms}

Lemma~\ref{lem:remove} and Proposition~\ref{prop:LL}
lead us to associate solutions to equation~(\ref{eqn:grn})
with $p$ satisfying condition~(\ref{eqn:cond}),
with newforms of certain levels. If there are
no newforms of the predicted levels, we immediately deduce that there
are no solutions to equation~(\ref{eqn:grn}).
With the help of a {\tt MAGMA} program we found all values of $D$ in
the range $1 \leq D \leq 100$ where there are no newforms
at the predicted levels. We deduce the following result.
\begin{cor}
Let $D$ be an integer belonging to the set
\[
	4,16,32,36,64.
\]
Then the equation~(\ref{eqn:grn}) does not have any solutions with
$p$ satisfying condition~(\ref{eqn:cond}).
\end{cor}
This Corollary does not add anything
new, since equation~(\ref{eqn:first}) has already been solved
by Cohn's method for $D=4,16,32,36, 64$ 
(but see \cite{Iv}, \cite{Siksek}, \cite{Le}). 
\end{ex}

\begin{ex}
Corollary~\ref{cor:Cohn} solves equation~(\ref{eqn:p}) 
for all values of $D$ in the range~(\ref{eqn:Drange})
except for $21$ values; these are the $19$ values listed
in~(\ref{eqn:left}) plus $D=55,~95$. As indicated
in Section~\ref{sec:survey} the cases $D=55,~95$ 
have been solved by Bennett and Skinner. It is however
helpful to look at the case $D=95$ again as it shows
how Methods I,~III work together in harmony. 
There is only one possible
signature $(d_1,d_2)=(1,95)$. Thus $t=x$, $s=y$
and we need to solve the equation
\begin{equation}\label{eqn:95}
t^2+95=s^p, 
\end{equation}
under the assumption that $p \geq 7$. 

Since $d_1=1$, it follows from Corollary~\ref{cor:Cohn}
that $y$ is even, and so $t=x$ is odd. Replacing $t$ by $-t$ if necessary,
we can assume that $t \equiv 1 \pmod{4}$. Table~\ref{table:LL1} leads us
to associate the solution $(t,s,p)$ with the Frey curve
\[
E_t~: \quad Y^2 + XY = 
	X^3 + \left(\frac{t - 1}{4}\right) X^2 + \left( \frac{t^2 + 95}{64} \right)X.
\] 
>From Proposition~\ref{prop:LL}, we know that any solution to equation~(\ref{eqn:95}) 
arises from a newform of level $190$. Using {\tt MAGMA} we find that there are,
upto Galois conjugacy, precisely four newforms at level $190$. These are
\begin{flalign*}
f_1 &= q - q^2 - q^3 + q^4 - q^5 + q^6 - q^7 + O(q^8), \\
f_2 &= q + q^2 - 3\, q^3 + q^4 - q^5 - 3\, q^6 - 5\, q^7 + O(q^8), \\
f_3 &= q + q^2 + q^3 + q^4 + q^5 + q^6 - q^7 + O(q^8),\\
f_4 &= q - q^2 + \phi\, q^3 + q^4 + q^5 - \phi\, q^6 + \phi\, q^7 + O(q^8), 
	\quad \text{where $\phi^2+\phi-4=0 $}.
\end{flalign*}
The first three newforms above are rational, and so correspond to the
three isogeny classes of elliptic curves of conductor $190$. It turns out
that none of these elliptic curves have non-trivial $2$-torsion.
By the remarks made after Proposition~\ref{prop:methodI} we know
that Method I will be successful in eliminating all but finitely many exponents $p$.
Indeed we find (in the notation of Proposition~\ref{prop:methodI}) 
that $B_3(f_1)=B_3(f_3)=15$. Thus we know that no solutions to
equation~(\ref{eqn:95}) arise from the newforms $f_1$ or $f_3$, since
otherwise, by Proposition~\ref{prop:methodI}, $p \mid 15$ which
contradicts our assumption that $p \geq 7$. We also find that 
$B_3(f_4)=2^4 \times 3$ and $B_7(f_4)=2^4 \times 7$. Thus
no solution arises from $f_4$. However, 
\begin{align*}
B_3(f_2) & =3 \times 7, & B_7(f_2) & =3^2 \times 5 \times 7, & B_{11}(f_2) &=0, \\
	B_{13}(f_2) & = 3 \times 5 \times 7 \times 13,  
	& B_{17}(f_2) & =3^2 \times 7 \times 11.
\end{align*}
We deduce that there are no solutions arising from $f_2$ with exponent
$p >7$. It does however seem likely that there is a solution with $p=7$.
Moreover, an attempt to prove that there is no solution with $p=7$
using Method~II fails: we did not find any integer $2 \leq n \leq 100$ 
satisfying the conditions of~\ref{prop:methodII}.

We apply Method III (and follow the notation of Section~\ref{sec:methodIII}). 
Write $\omega=\frac{1+\sqrt{-95}}{2}$.
Taking 
\[
S=\left\{ 113, 127, 239, 337, 491  \right\}
\]
we find that 
\[
\cap_{l \in S} \Gamma_l= \left\{ \frac{-528 -2 \omega}{2187}  \right\}.
\]
Thus if we have any solutions at all then, by Proposition~\ref{prop:methodIII1},
we know
\[
(t-1) + 2 \omega = \left( \frac{-528 -2 \omega}{2187}\right) 
	\left(U + V \omega \right)^7,
\] 
for some integers $U,~V$. Equating imaginary parts and clearing the denominators
we find that  
\[
\begin{split}
- U^7 & - 1855 V U^6 - 5061 V^2 U^5 + 214165 V^3 U^4 + 416605 V^4 U^3 \\
	& - 2834013 V^5 U^2 - 2944375 V^6 U + 2818247 V^7= 2187.
\end{split}
\]
Using {\tt pari/gp} we find that the only solution to this Thue equation is given by %%!
$U=-3$, $V=0$. %%!
This shows that $(t,s)=(529,6)$. 

The reader will notice that $(t,s)=(-529,6)$ is also a solution
to equation~(\ref{eqn:95}) with $p=7$ but it seems to have been 
\lq missed\rq\ by
the method. This is not the case; we are assuming that the sign of 
$t$ has been chosen so that $t \equiv 1 \pmod{4}$. The solution
$(t,s)=(-529,6)$ arises from some other newform (probably at some
different level) and via a different Frey curve which we have not
determined.
\end{ex}

\begin{ex}
For our last example we look at the case where $D=25$.
This, like $18$ other cases, must be resolved by a combination
of the modular approach and our lower bound for linear forms
in three logarithms which is to come.
We assume that $p \geq 7$, and so $p$ satisfies conditions~(\ref{eqn:cond}).
There are now two possible signatures 
$(d_1,d_2)=(1,25),~(5,1)$ satisfying the conditions of 
Lemma~\ref{lem:remove}.
However, by Corollary~\ref{cor:Cohn}, we may suppose that 
$d_1 > 1$ and so $d_1=5$, $d_2=1$. %%!!
We write $t=x/5$, $s=y/5$ where we know that $t,~s$ are
integral by Lemma~\ref{lem:remove}. Equation~(\ref{eqn:st1})
becomes
\[
t^2+1=5^{p-2} s^p, \quad t \neq \pm 1.
\]
Following Table~\ref{table:LL1}, we associate with any solution to this
equation the Frey curve
\[
E_t: \qquad Y^2=X^3+2t X^2 -X,
\]
and we know by Proposition~\ref{prop:LL} that any solution
must arise from a newform of level $160$. 
Using the computer algebra system {\tt MAGMA}
we  find that there are upto Galois conjugacy three
such newforms:  
\begin{flalign*}
f_1 &= q - 2 q^3 - q^5 - 2 q^7 + O(q^8), \\
f_2 &= q + 2 q^3 - q^5 + 2 q^7 + O(q^8), \\
f_3 &= q+2 \sqrt{2} q^3 + q^5 - 2 \sqrt{2} q^7 + O(q^8).
\end{flalign*}
The first two newforms are rational, corresponding
respectively to elliptic curves 160A1 and 160B1
in Cremona's tables~\cite{Cre}. The third has 
coefficients in $K=\Q(\sqrt{2})$ and  is 
straightforward to eliminate using Method I. In the notation of 
Proposition~\ref{prop:methodI} we find that if $f_3$ does give
rise to any solutions $(t,s,p)$ then $p \mid B_3(f_3)=24$.
This is impossible as $p \geq 7$, and so $f_3$ does not give 
rise to any solutions.

We where unable to eliminate newforms $f_1,~f_2$ using Method I.
Instead using our implementation of Method II in {\tt MAGMA}
we showed that there are no solutions arising from either form
with $7 \leq p \leq 100$.  With our implementation of the improved 
Method II (Proposition~\ref{prop:methodIIrational})
in {\tt pari/gp} we showed that there are no solutions with 
$100 \leq p \leq 163762845$; this took roughly 26 hours on 2.4 GHz
Pentium IV PC. The choice of where to stop the computation is of course
not arbitrary, but comes out of our bound for the linear form in logarithms.
We will later prove that $p \leq 163762845$ thereby completing the resolution
of this case.
\end{ex}

\section{Results III}\label{sec:results}

We applied the methods of the previous sections to
solve all equations~(\ref{eqn:p})
with $D$ is our range~(\ref{eqn:Drange}).

\begin{table}
\begin{minipage}{\linewidth}
\caption{Computational details for Lemma~\ref{lem:remaining} and its proof.}
\begin{tabular}{|c|c|l|c|c|l|}
\hline
$D$ & $(d_1,d_2)$ & $E$~\footnote{We give here the Cremona code for the elliptic curves $E$ as 
in his book \cite{Cre} and his online tables:\ \  
{\tt http://www.maths.nott.ac.uk/personal/jec/ftp/data/INDEX.html}} & $p_0$ & 

Machine~\footnote{The machines are as follows
	\begin{enumerate} 
		\item[{\bf P1}] 2.2 GHz Intel Pentium PC.
		\item[{\bf P2}] 2.4 GHz Intel Pentium PC.
		\item[{\bf S1}] Dual processor 750MHz UltraSPARC III.  
		\item[{\bf S2}] 650 MHz UltraSPARC IIe.
		\item[{\bf S3}] UltraSPARCIII with 12 processors of 1050 MHz speed.
	\end{enumerate}
} 

& Time\\
\hline\hline
$7$ & $(1, 7)$ & 14A1 & $181\, 000\, 000$ & P1  & 26h, 43mn \\
\hline
$15$ & $(1, 15)$ & 30A1 & $624\, 271\, 465$ & S1 & 252h, 50mn \\
\hline
$18$ & $(3, 2)$ & 384D1, 384A1, & $306\, 111\, 726$     &S3     & 293h, 14mn\\
        &       & 384G1, 384H1  &       &       & \\
\hline
$23$ & $(1, 23)$ & 46A1 & $855\, 632\, 066$ &  S2 &  477h, 36mn\\
\hline
$25$ & $(5, 1)$ & 160A1, 160B1 & $163\, 762\, 845$ & P2 & 25h, 58mn \\
\hline
$28$ & $(2, 7)$  & 14A1 & $315\, 277\, 186$     & P1    & 55h, 41mn \\
\hline
$31$ & $(1, 31)$ &  62A1 & $860\, 111\, 230$    & S3    & 242h, 2mn \\
\hline
$39$ & $(1, 39)$ & 78A1 & $852\, 830\, 725$ & P1  & 193h, 41mn\\
\hline
$45$ & $(3, 5)$ & 480B1, 480F1, & $340\, 749\, 424$ & S1  & 448h, 43mn\\
        &       & 480G1, 480H1  &       &       &\\
\hline
$47$ & $(1, 47)$ & 94A1 & $1\, 555\, 437\, 629$ & S3      & 451h, 34mn\\
\hline
$60$ & $(2, 15)$ & 30A1 & $358\, 541\, 296$ & S1 &  130h, 30mn\\
\hline
$63$ & $(1, 63)$ & 42A1 & $292\, 825\, 735 $ & S1 & 99h, 45mn  \\
\hline
$71$ & $(1, 71)$ & 142C1 & $2\, 343\, 468\, 548$ & S3     & 697h, 26mn \\
\hline
$72$ & $(3, 8)$ & 96A1, 96B1 & $451\, 620\, 034$ & S1   & 316h, 27mn \\
\hline
$79$ & $(1, 79)$ & 158E1 & $1\, 544 \, 381\, 661$       & S3      & 448h, 47mn      \\
\hline
$87$ & $(1, 87)$ & 174D1 & $1\, 148\, 842 \, 108 $      & S3    & 329h, 45mn    \\
\hline
$92$ & $(2, 23)$ & 46A1 & $996\, 255\, 151$&    S3 & 285h, 10mn       \\
\hline
$99$ & $(3, 11)$ & 1056B1, 1056F1 & $593\,734\,622$ & P2 & 138h, 46mn   \\
\hline
$100$ & $(5, 4)$ & 20A1 & $163\, 762\, 845$ & P1 & 21h, 23mn    \\
\hline
\end{tabular}
\label{table:rational}
\end{minipage}
\end{table}

\begin{lem}\label{lem:remaining}
Suppose $D$ is in the range~(\ref{eqn:Drange}) and $p$ is a prime
satisfying~(\ref{eqn:cond}). Suppose $(x,y,p)$ is a solution to 
equation~(\ref{eqn:p}) that is not included in the tables.
Then $D$ is one of 
\begin{equation}\label{eqn:remaining}
7, 15, 18, 23, 25, 28, 31, 39, 45, 47, 60, 63, 71, 72, 79, 87, 92, 99, 100.
\end{equation}
Moreover $(x,y,p)$ has signature $(d_1,d_2)$ and arises from an elliptic curve 
$E$ and $p > p_0$
where $E$, $p_0$ and $(d_1,d_2)$ are given by Table~\ref{table:rational}.
\end{lem}
\begin{proof}
We wrote a {\tt MAGMA} program that does the following:
For each $D$ in the range~(\ref{eqn:Drange}) we write down 
the set of possible signatures $(d_1,d_2)$ satisfying
the conditions of Lemma~\ref{lem:remove}.

For each such pair $(d_1,d_2)$ write down the (one or two) Frey curves 
given by the Tables~\ref{table:LL1}, \ref{table:LL2}, \ref{table:LL3},
bearing in mind the information given by Corollary~\ref{cor:Cohn}.

For each Frey curve we compute the conductor (given by Proposition~\ref{prop:LL}) 
of the newforms giving rise to possible solutions,
and then write down all these newforms. 

We attempt to eliminate each newform $f$ using Method I. This involves searching
for primes $l \nmid 2D$ such that 
(in the notation of Proposition~\ref{prop:methodI}) $B_l(f) \neq 0$.
If we are successful and find such primes $l_1,\ldots, l_m$ then
by Proposition~\ref{prop:methodI} the exponent $p$ divides 
all the $B_{l_i}(f)$, and so divides their greatest common divisor $B$ (say).
If $B$ is divisible by any prime $p$ that satisfies condition~(\ref{eqn:cond})
then we attempt to eliminate this possible exponent $p$ 
using Method II: this involves searching for an integer $2 \leq n \leq 100$ satisfying 
conditions (a), (b), (c) of Proposition~\ref{prop:methodII}. If one such $n$ is
found then we know that there are no solutions for the particular exponent $p$.
Otherwise we apply Method III 
(Proposition~\ref{prop:methodIII1})
 to write down Thue equations leading to possible solutions (see below). 

%In a handful
%of cases we used the alternative method of Section~\ref{sec:alter} to show
%that the equation $x^2+D=y^p$ has no solutions.

As predicted by the comments made after Proposition~\ref{prop:methodI},
Method I succeeded with all non-rational newforms and all rational newforms
corresponding to elliptic curves with only trivial $2$-torsion (it also
succeeded with some rational newforms corresponding to elliptic curves with
non-trivial $2$-torsion). Indeed, we found no solutions arising from non-rational
newforms for $D$ in our range~\ref{eqn:Drange}.

We are left only with rational newforms $f$ that correspond to elliptic curves $E$
having some non-trivial $2$-torsion. The details of these are documented
in Table~\ref{table:rational}. For primes 
$p < 100$ satisfying condition~(\ref{eqn:cond}) we attempt to
show that there are no solutions arising from $E$ for the particular
exponent $p$ using Method II (as before). If this fails for a particular
exponent $p$, then we use Method III to write down the Thue equations leading
to the possible solutions.

Our proof that $p \geq 100$ is now complete except that  
there are some Thue equations to solve.
We had to solve
Thue equations of degree $7$ for $D=7$, $47$, $79$ and $95$. %%!!!
These were solved
using {\tt pari/gp} and the solutions incorporated in our Tables. We also had
to solve a Thue equation of degree $11$ for $D=23$, of degree $17$ for $D=28$,
and of degree $13$ for $D=92$. We were unable to (unconditionally) solve these
three Thue equations using the in-built functions of {\tt pari/gp}. The
reason is that, in each case, it was impossible for {\tt pari/gp} to
prove that the system of units it had found -- though of correct rank -- was maximal.
We are grateful
to Dr. Guillaume Hanrot for sending us his {\tt pari} program for
solving Thue equations without the full unit group.
This program, based on~\cite{Hanrot}, solved all three equations in a few
minutes.

%Again, for a handful
%of cases we applied the method explained in Section~\ref{sec:alter} to show
%that the equation $x^2+D=y^p$ has no solutions. 

For the next step we implemented our improved Method II
(Proposition~\ref{prop:methodIIrational}) in 
{\tt pari/gp} (see the remark after the proof).
To complete the task and show that $p > p_0$ for any missing solution
we used our {\tt pari/gp} program to disprove
the existence of any missing solution for each prime $100 \leq p \leq p_0$. 
We ran this {\tt pari/gp} program on various
machines as indicated in Table~\ref{table:rational}. The total computer time for
this step is roughly $206$ days.
\end{proof}

\begin{remark}
The reader may be surprized that some of the computations
were done in {\tt MAGMA} while others were carried out
in {\tt pari/gp}. As stated earlier, {\tt MAGMA} has a package
for computing modular forms. This is essential for us, and is
unavailable in {\tt pari/gp}.

For showing that $p > p_0$, it is simply not practical to
use {\tt MAGMA}. Here we are using the improved Method~II
(Proposition~\ref{prop:methodIIrational}). 
The main bottle-neck in Method II is computing $a_l(E)$
for primes $l$ that can be about $10^{11}$ (recall $l$ is
a prime satisfying $l \equiv 1 \pmod{p}$). For this {\tt pari/gp}
uses the theoretically slower Shanks-Mestre method \cite{Cohen}
rather than the theoretically faster Schoof-Elkies-Atkin \cite{Schoof}  
method used by {\tt MAGMA}.
But for primes of the indicated size it seems that {\tt pari/gp}
is about 10 times faster than {\tt MAGMA}.

The reader may also note that two of the machines we used are multiprocessor
machines. The computation for each $D$ could have been speeded up considerably by parallelising.
We however decided against this, so as to keep our programs simple and transparent.
\end{remark}

\section{The \lq Modular\rq\ Lower Bound for $y$}

In this section we would like to use the modular approach to prove
a lower bound for $y$ with $D$ in the range~(\ref{eqn:Drange}).
Before doing this we prove a general result for arbitrary non-zero
$D$. 

\begin{prop}\label{prop:LB}
Suppose $D$ is a non-zero integer, and $d_1$, $d_2$ satisfy (i)--(v) %%!
of Lemma~\ref{lem:remove}. Suppose $(t,s,p)$ is a solution
to equation~(\ref{eqn:st1}) arising from a {\bf rational} newform
$f$ via a Frey curve $E_t$. Then either $\mathrm{rad} (s) \mid 2 d_1$ 
%% changed {\rm rad} to \mathrm{rad}
or $\lvert s \rvert > (\sqrt{p} - 1)^2$.
\end{prop}
\begin{proof}
Since the newform is rational we know that the newform $f$
corresponds to an elliptic curve $E/\Q$ whose conductor
equals the level of $f$.

Suppose $\mathrm{rad}(s)$ does not divide $2d_1$. Since $t$ %% again
and $d_2$ are coprime we see that there is some prime $l \mid s$
so that $l \nmid 2D$. By Lemma~\ref{lem:cong2} we see that 
$p$ divides $l+1  \pm a_l(E)$. It follows from the Hasse--Weil 
bound that $l+1  \pm a_l(E) \neq 0$, and so
\[
p \leq l+1 \pm a_l(E) < (\sqrt{l}+1)^2,
\]
where again we have used Hasse--Weil. Thus $l > (\sqrt{p} - 1)^2$. 
The Proposition follows as $l\mid s$.  %%!
\end{proof}

\begin{cor}\label{cor:LB}
Suppose $D$ is one of the values in~(\ref{eqn:remaining}). %%! 
If $(x,y,p)$ is a solution to equation~(\ref{eqn:grn})
not in the Tables below then $y > (\sqrt{p} - 1)^2$.
\end{cor}
\begin{proof}
Suppose $D$ is in the range~(\ref{eqn:Drange})
and $(x,y,p)$ is some solution to equation~(\ref{eqn:grn})
not in the below Tables. From the preceding sections
we know that this solution must satisfy condition~(\ref{eqn:cond}).
Moreover by Lemma~\ref{lem:remove}, 
\[
	x=d_1 t,\qquad y=\mathrm{rad}(d_1) s, %% changed {\rm rad} to \mathrm{rad}
\] 
where $(t,s,p)$ satisfy equation~(\ref{eqn:st1})
for some $d_1,~d_2$ satisfying conditions (i)--(v) of that Lemma.

We have determined for $D$ in the specified range all solutions
to equations~(\ref{eqn:st1}) arising from non-rational newforms (indeed there were none).
Thus we may suppose that our putative %phantom 
solution arises from 
a rational newform. By Proposition~\ref{prop:LB} 
we see that either $\lvert y \rvert \geq \lvert s \rvert > (\sqrt{p} - 1)^2$
or $\mathrm{rad} (s) \mid 2 d_1$. We must prove that the second possibility
does not arise.

Suppose that $\mathrm{rad} (s) \mid 2 d_1$. From Lemma~\ref{lem:remove}
we see that $\mathrm{rad}(y) \mid 2 d_1$. We first show that 
$\mathrm{rad}(y) \neq 2$. For in this case we have reduced to 
an equation of the form $x^2+D=2^m$. For $\lvert D \rvert < 2^{96}$,
which is certainly the case in our situation, Beukers \cite[Corollary 2]{Be}
shows that 
\[
m \leq 18+2 \log{\lvert D \rvert} / \log{2}.
\]
A short {\tt MAGMA} program
leads us to all the solutions to this equation for $2 \leq D \leq 100$
and we find that these are already in our tables.

Thus we may suppose that $\mathrm{rad}(y) \mid 2 d_1$ and 
$\mathrm{rad}(y) \neq 2$. An examination of the possible cases reveals
the following possibilities
\begin{enumerate}
\item $D=18$, $45$, $72$, $99$ and $\mathrm{rad}(y)=3$,
\item $D=25$, $100$ and $\mathrm{rad}(y)=5$.
\end{enumerate}
On removing the common factors, each case quickly reduces
to an equation that has already been solved. For
example, we must solve $x^2+100=y^p$ under
the assumption that $\mathrm{rad}(y)=5$ or equivalently
 the equation $x^2+100=5^m$. Removing the common
factor reduces to the equation $X^2+4=5^{m-2}$.
But the equation $X^2+4=Y^n$ has already been solved
and has only the solutions $(X,Y,n)=(2,2,3)$, $(11,5,3)$.
We quickly see that the only solution
to $x^2+100=y^p$ with $\mathrm{y}=5$ is $(x,y,p)=(55,5,5)$.
\end{proof}

%\input amssym.def
%\input amssym.tex

%\let\eps=\varepsilon

%\def\set#1{\left\{#1\right\}}

%\def\cf{{\it cf}}
%\def\eg{{\it e.g}}
%\def\ie{{\it i.e}}
 
%\def\tmod{\kern -4 pt \pmod}
%\def\div{\mathbin{\vert}}
%\def\order{{\rm order}\kern .7 pt}
%\def\notdiv{\mathop{\not \hskip 0.9 mm\mid}\nolimits}
%\def\prf{\noindent{\sl Proof. ---\ }}
%\def\kro#1#2{\left({#1 \over #2}\right)}

%\def\qed{\lower 2pt\hbox{\vrule\vbox to 10pt{\hrule width 4pt
%                                             \vfill
%                                             \hrule}\vrule}}

%\def\cqfd{\unskip\penalty 500\kern 10pt\qed}
%\long\def\fo#1\of{}

%\def\bR{{\Bbb R}}
%\def\bC{{\Bbb C}}
%\def\bF{{\Bbb F}}
%\def\bK{{\Bbb K}}
%\def\bQ{{\Bbb Q}}
%\def\bZ{{\Bbb Z}}
%\def\bN{{\Bbb N}}  

%\def\cA{{\cal A}}
%\def\cI{{\cal I}}
%\def\cM{{\cal M}}
%\def\cN{{\cal N}}
%\def\cS{{\cal S}}

%\def\h{\mathop{\rm h\kern 1pt}}
%\def\H{\mathop{\rm H\kern 1pt}}
%\def\L{\mathop{\rm L\kern 1pt}}
%\def\M{\mathop{\rm M\kern 1pt}\nolimits}
%\def\N{\mathop{\rm N\kern 1pt}}
%\def\v{\mathop{\rm v}}
%\def\w{\mathop{\rm w}}
%\def\card{{\rm Card}\kern.7pt}
\def\Log{{\rm Log}\kern.7pt}
%\def\sign{{\rm sign}\kern.7pt}

%\def\gA{{\frak A}}
%\def\ga{{\frak a}}
%\def\gb{{\frak b}}
%\def\gc{{\frak c}}
%\def\gj{{\frak j}}
%\def\gp{{\frak p}}
%\def\gq{{\frak q}}
%\def\gd{{\frak d}}
%\def\gS{{\frak S}}

%\long\def\fo#1\of{}

%%%%%%%%%%%%%%%%

%\let\eps=\varepsilon

%\def\set#1{\left\{#1\right\}}

%\def\cf{{\it cf}}
%\def\eg{{\it e.g}}
%\def\ie{{\it i.e}}
 
%\def\tmod{\kern -4 pt \pmod}
%\def\div{\mathbin{\vert}}
%\def\order{{\rm order}\kern .7 pt}
%\def\notdiv{\mathop{\not \hskip 0.9 mm\mid}\nolimits}
%\def\prf{\noindent{\sl Proof. ---\ }}
%\def\kro#1#2{\left({#1 \over #2}\right)}

%\def\qed{\lower 2pt\hbox{\vrule\vbox to 10pt{\hrule width 4pt
%                                             \vfill
%                                             \hrule}\vrule}}

%\def\cqfd{\unskip\penalty 500\kern 10pt\qed}
%\long\def\fo#1\of{}

%\def\bR{{\Bbb R}}
%\def\bC{{\Bbb C}}
%\def\bF{{\Bbb F}}
%\def\bK{{\Bbb K}}
%\def\bQ{{\Bbb Q}}
%\def\bZ{{\Bbb Z}}
%\def\bN{{\Bbb N}}  

%\def\cA{{\cal A}}
%\def\cI{{\cal I}}
%\def\cM{{\cal M}}
%\def\cN{{\cal N}}
%\def\cS{{\cal S}}

%\def\h{\mathop{\rm h\kern 1pt}}
%\def\H{\mathop{\rm H\kern 1pt}}
%\def\L{\mathop{\rm L\kern 1pt}}
%\def\M{\mathop{\rm M\kern 1pt}\nolimits}
%\def\N{\mathop{\rm N\kern 1pt}}
%\def\v{{\rm v}}
%\def\w{\mathop{\rm w}}

%\def\gA{{\frak A}}
%\def\ga{{\frak a}}
%\def\gb{{\frak b}}
%\def\gc{{\frak c}}
%\def\gd{{\frak d}}
%\def\gg{{\frak g}}
%\def\gj{{\frak j}}
%\def\gp{{\frak p}}
%\def\gq{{\frak q}}
%\def\gr{{\frak r}}
%\def\gS{{\frak S}}

%\long\def\fo#1\of{}

%\magnification 960
%\hsize 172mm
%\vsize 210mm

%\vglue -11mm

%\vglue 2cm

%\centerline{\bf On the equation $x^2+C=y^p$}
 
%\vglue 2cm

%\noindent
\section{The Linear Form in Logarithms}~\label{sec:lin}

It is useful at this point to recap what we have done so far.
We would like to complete the proof of Theorem~\ref{thm:main} by showing that
our Tables at the end are not missing any solutions.
So let us suppose that our Tables at the end are missing some solution
$(x,y,p)$ to equation~(\ref{eqn:p}) for some value of $D$ in our range~(\ref{eqn:Drange}).
We have proved (Lemma~\ref{lem:remaining}) 
that $D$ is one of the values in~(\ref{eqn:remaining}).
Moreover, (again by Lemma~\ref{lem:remaining} and by Corollary~\ref{cor:LB}) 
any missing solution $(x,y,p)$ must satisfy 
\begin{equation}\label{eqn:ylb}
p > p_0, \qquad  y \geq (\sqrt{p}-1)^2,
\end{equation}
with $p_0$ given by Table~\ref{table:rational}. 
Our aim is to show that $p \leq p_0$ thus obtaining a contradiction.

>From the table of values of $p_0$ we know that
\begin{equation}\label{eqn:xypLB}
\lvert x \rvert ,~y,~p \geq 10^8
\end{equation}
and indeed much more, though this inequality is sufficient for much of our
later work. In the remainder of this paper we assume that $D$ is one of the
remaining values~(\ref{eqn:remaining}), and always write (as before)
$D=D_1^2 D_2$, where $D_2$ is square-free. The triple $(x,y,p)$ will always
be a solution to equation~(\ref{eqn:p}) supposedly missing from our Tables
and hence satisfying the above inequalities.

In this section we write down the linear form in logarithms corresponding to
the equation~(\ref{eqn:p}) and apply a Theorem of Matveev to obtain upper bounds for 
the exponent $p$. These upper bounds obtained from Matveev's Theorem are not small enough
to contradict our lower bounds for $p$ obtained in Lemma~\ref{lem:remaining}
but they are needed when we come to apply our bounds for linear forms in three logarithms given
in the next section.

%\begin{table}
%\caption{}
%\begin{tabular}{||c|c|c||c|c|c||c|c|c||}
%\hline\hline
%$D$ & $d$ & $k_0$ & $D$ & $d$ & $k_0$ & $D$ & $d$ & $k_0$\\
%\hline\hline
%$7$ &  $2$ & $1$ & $31$ & $2$ & &$71$ & $2$ & \\
%\hline
%$15$ & $2$ & $1$ & $39$ & $2$ & &$72$ & $3$ & \\
%\hline
%$18$ & $3$ & $3$ & $45$ & $3$ & &$79$ & $2$ & \\
%\hline
%$23$ & $2$ & $3$ & $47$ & $2$ & &$87$ & $2$ & \\
%\hline
%$25$ & $5$ & 	& $60$  &  $4$ & &$92$ & $4$ & \\
%\hline
%$28$ &  $4$ & $3$ & $63$ & $2$ & &$99$ & $3$ & \\
%\hline
%$100$ & $5$ &     &       &     & &       &  &  \\
%\hline\hline
%\end{tabular}
%\label{table:kappa}
%\end{table}
%

\begin{lem}\label{lem:kappa}
Let $(d_1,d_2)$ be the signature of our supposedly missing solution $(x,y,p)$
(which we know from Lemma~\ref{lem:remaining}).
Define
\begin{equation}\label{eqn:d}
d=
\begin{cases}
d_1, & \text{for $D \not \equiv 7 \pmod{8}$,} \\
2 d_1, & \text{for $D  \equiv 7 \pmod{8}$.}
\end{cases}
\end{equation}
Then $d$ is a prime power, say $d=q^c$ for prime $q$,
where moreover, $q$ splits in $\cL=\Q(\sqrt{-D_2})$,
say $(q)=\gq \bar\gq$.
Let $k_0$ be the smallest positive integer such
that the ideal $\bar\gq^{ k_0}$ is principal, %%!
say $\bar\gq^{ k_0}=(\alpha_0) $. %%!
Also let
\[
k=\begin{cases}
k_0, & \text{if $k_0$ is odd}, \\
k_0/2, & \text{if $k_0$ is even}, %%!
\end{cases}
\quad \text{and}\quad
\kappa=  \begin{cases}
2, & \text{if $k_0$ is odd,}\\
1, & \text{if $k_0$ is even,}
\end{cases}
\quad \text{so that\ } k=\frac{\kappa k_0}{2}.
\]
Then there exists $\gamma \in \cL$ such that 
\[
	\left( \frac{x-D_1 \sqrt{-D_2}}{x+D_1\sqrt{-D_2}}\right)^{k}
	=\alpha^\kappa \gamma^p,
\]
where
\[	
	\alpha=\bar\alpha_0 /\alpha_0, \quad
	\h(\alpha)=\frac{k\log d}{\kappa}, \quad \h(\gamma)= \frac{k\log y}{2}.
\]
\end{lem} 
\begin{proof}
% Let $h_\bK=h$ be the class-number of this field; we suppose $n$ coprime with $h$.
We begin with the factorization 
\[
(x+D_1\sqrt{-D_2})(x-D_1\sqrt{-D_2})=y^p.
\]
Our first step is to show that any prime divisor $q$ of $y$ splits in $\cL$.
Suppose otherwise, then we may write $(q)=\gq$ or $(q)=\gq^2$ for some prime
ideal $\gq$ satisfying $\bar \gq = \gq$. If $p=2r+1$ then clearly
$\gq^r$ divides both factors on the left-hand side above, and so
divides  $2D_1 \sqrt{-D_2}$. This is impossible in view of the fact
that $p$ is enormous, and $2 \leq D \leq 100$. 
Thus we have shown that every prime divisor
$q$ of $y$ splits in $\cL$.
 
Let us write
\[
	y = \prod_{i\in  I} {q_i}^{a_i}\qquad
	\mathrm{and}\qquad (q_i)=\gq_i \bar \gq_i, \quad \gq_i \not=\bar \gq_i,\ i\in I.
\]
Then
\[
	(x+D_1 \sqrt{-D_2}) = \prod_{i\in  I} ({\gq_i}^{b_i} {\bar \gq_i}^{c_i}),
\]
where we assume (for commodity of notation) that $b_i\ge c_i$ for all $i$, and thus
\[
	(x-D_1 \sqrt{-D_2}) = \prod_{i\in  I}({\gq_i}^{c_i} {\bar \gq_i}^{b_i}),
\]
with
\[
	b_i+c_i = p a_i, \quad \hbox{for all $\,i\in I$.}
\]

Let
\[
	\gd = \gcd \left( x+D_1 \sqrt{-D_2}, x-D_1 \sqrt{-D_2} \right)\, ;
\]
clearly
\[
	\gd = \prod_{i \in I} ({\gq_i} {\bar \gq_i})^{c_i}= \prod_{i\in  I} (q_i)^{c_i} . 
\]
%
%Moreover, since $\gd$
%divides $y^p$ and all the prime divisors of $y$ split in $\cL$,
%we see that $\gcd(\gd,D_\cL)=\gcd(\gd,D_2)=1$
%where $D_\cL$ is the discriminant of the field~$\cL$.
%
%
This shows that $\gd=(d)$ where $d\in \Z$. We would like to calculate
this $d$ and verify that its value is in agreement with~(\ref{eqn:d}).
>From the definition of $\gd$ we see that $d \mid 2x$ and $d \mid 2 D_1$.
However, by our definition of signature, $\gcd (x^2,D)=d_1^2$. 
It follows that $d^2 \mid 4 d_1^2$ and so
$d \mid 2 d_1$. But $d_1 \mid x$ and $d_1 \mid D_1$. Hence $d_1 \mid \gd$
and so $d_1 \mid d$. Thus $d=d_1$ or $d=2d_1$.

%We recall by Corollary~\ref{cor:Cohn} that either $D \equiv 7 \pmod{8}$
%and $2 \mid y$ or $d_1 > 1$. 
We note the following cases:
\begin{itemize}
\item If $D_2 \not \equiv 7 \pmod{8}$ then $2 \nmid y$. 
	Thus $2 \nmid d$ and so $d=d_1$.
\item Suppose $D_2 \equiv 7 \pmod{8}$. 
	Now from Lemma~\ref{lem:remove} and its proof we know that $D=d_1^2 d_2$ 
	and $x=d_1 t$ where $\gcd (t,d_2)=\gcd(d_1,d_2)=1$. 
	Clearly $d_2=d_3^2 D_2$ with $d_3=D_1/d_1$ integral. 
Suppose first that $d_1$ is even. It follows easily that $t,d_2$
	are odd and
	\[
		(d)=\gd=  
	% d_1(t+d_3\sqrt{-D_2},t-d_3\sqrt{-D_2})=
		2 d_1 \Biggl( \frac{t+d_3\sqrt{-D_2}}{2},\frac{t-d_3\sqrt{-D_2}}{2} \Biggr).
	\]
	Hence $(2 d_1) \mid d$ and so $d=2 d_1$.
\item The only case left to consider is $D_2 \equiv 7 \pmod{8}$ and $d_1$ is odd.
	By examining Table~\ref{table:rational} we see that $d_1=1$. Thus $2 \mid y$ 
	by Corollary~\ref{cor:Cohn}. Clearly $x$ is odd, and 
	the same argument as above shows that $d=2=2 d_1$.
\end{itemize}
This proves that $d$ satisfies~(\ref{eqn:d}). By looking again at the possible
values of $d_1$ in Table~\ref{table:rational} we
see that $d$ is a prime-power in all cases. Let $j \in I$ such that
$d=q_j^{c_j}$. Thus $c_i=0$ for all $j \neq i$. Then
\[
	(x+D_1 \sqrt{-D_2}) = {\bar\gq_j}^{c_j} 
	\cdot {\gq_j}^{b_j} \cdot \prod_{j \neq i}{\gq_i}^{p a_i},
\]
whence
\[
	(x+D_1 \sqrt{-D_2}) = (\bar\gq_j\,{\gq_j}^{-1})^{c_j} 
	\cdot \prod_{i\in  I} {\gq_i}^{pa_i}
	=(\ga\,{\bar\ga}^{-1})\, \gg^p,
\]
where $\ga$ and $\gg$ are integral ideals with
\[
	\ga= {\bar\gq_j}^{c_j},
	\quad \N (\ga) = {q_j}^{c_j} = d , 
	\quad \N (\gg )= y,
\]
and $\N$ denotes the norm. Thus, as ideals, %% slight change
\[
	\left( \frac{x-D_1\sqrt{-D_2}}{x+D_1\sqrt{-D_2}}\right) = 
	(\bar\ga\,{\ga}^{-1})^2\, (\bar\gg\,{\gg}^{-1})^p.
\]

%We define $k_0$ to be  the smallest positive rational integer such 
%that the ideal $\ga^{k_0}$ is principal. In all cases, since we know $d$
%(and it is a prime power),
%we know $\ga$ upto conjugation, and we were able to find $k_0$ with the help
%of {\tt pari/gp}; the values of $k_0$ are tabulated in Table~\ref{table:kappa}.
We define $k_0,~k,~\kappa,~\alpha_0$ as in the statement of the Lemma.
Thus $\ga^{k_0} = (\alpha_0)$ and we have the relation (between ideals)
\[
	(x+D_1 \sqrt{-D_2})^{k}=(\ga/\bar \ga)^k \gg^{k p}=\ga^{2k}(\N(\ga))^{-k}\gg^{k p}
	=(\alpha_0)^\kappa (d)^{-k}\,\gg^{k p}.
\] 
However $p$ is an enormous prime certainly not dividing the class number. 
This shows that $\gg^{k}$ is also principal, 
$\gg^{k}=(\gamma_0)$,
say, where $\gamma_0$ is an algebraic integer chosen so that 
the following equality of elements of $\cL$ holds
\[
	(x+D_1 \sqrt{-D_2})^{k}=\alpha_0^\kappa d^{-k} \gamma_0^p.
\]
Note that
\[ 
	\N (\alpha_0)=d^{k_0}, \quad \N (\gamma_0)=y^{k}.
\]
Write
\[
	\alpha=\bar \alpha_0/ \alpha_0, \quad \gamma= \pm \bar\gamma_0/\gamma_0.
\]
The proof of the Lemma is complete except for
the statements about the heights of $\alpha,~\gamma$. These follow from Lemma~\ref{lem:M1}
below.
\end{proof}

\begin{lem}\label{lem:M1}
Let $\alpha$ be an algebraic number whose conjugates are all (including $\alpha$ itself) of
modulus equal to~1, then
\[
        \h (\alpha) =
	        \frac{1}{\deg \alpha}\log a,
\]
where $a$ is the leading coefficient of the minimal polynomial of $\alpha$.
In particular, if $\alpha=\bar\alpha_0/\alpha_0$ where $\alpha_0$ is a non-real quadratic
irrationality, then
\[
        \h (\alpha) ={\textstyle \frac{1}{2}}\log \N(\alpha_0).
\]
\end{lem}
\begin{proof}
Let $d=\deg \alpha$.
By hypothesis $\alpha$ is a root of a polynomial of $\Z[X]$ of the form
$P(X)=aX^d+\cdots$. We have
\[
        \h (\alpha)
       ={\textstyle \frac{1}{d}} \log  \M (P),
\]
where $\M$ denotes Mahler's measure, and the first result easily follows since the roots of $P$ are
of modulus equal to~1. This proves the first assertion.

The proof of the second assertion, which is quite easy, is omitted.
\end{proof}

%\noindent{\bf Example .} ---
%For $C=15$ we have $h=2$, then $k=1$ and from the above study it is easy to show
%that the Diophantine equation
%$$
%x^2+15=2^n
%$$
%has only the solution
%$$
%1^2+15=2^4.
%$$
We now write the linear form in three logarithms.
Define %% slight change
\[
\Lambda = \log \left( \frac{x-D_1 \sqrt{-D_2}}{x+D_1 \sqrt{-D_2}} \right),
\]
where we have taken the principal determination of the logarithm. 
\begin{lem}\label{lem:logLambda}
\[
	\log |\Lambda | \le - \frac{p}{2}\,\log y+ \log\bigl(2{.}2\,D_1 \sqrt {D_2} \bigr).
\]
\end{lem}
\begin{proof}
We will rely on the lower bounds~(\ref{eqn:xypLB}). Clearly
\[
	\left|\frac{x-D_1 \sqrt{-D_2}}{x+D_1 \sqrt{-D_2}}-1\right| < 2 \, 
	\frac{D_1 \sqrt{D_2}}{\lvert x \rvert}.
\]
A standard inequality (\cite{Smart}, Lemma B.2) shows that 
\[
	| \Lambda  | <  2{.}1 \, \frac{D_1 \sqrt {D_2}}{\lvert x \rvert},
\]
so that
\[
\log  | \Lambda  | <  -\log \lvert x \rvert + \log\bigl(2{.}1 \,D_1 \sqrt {|D_2|}\bigr).
\]
Using the fact that $y^p-x^2=D$, and a similar argument to the one above,
we deduce the Lemma.
\end{proof}

The main tool to bound $p$ will be the theory of linear forms of (at most three)
logarithms.  We need the special case of three logarithms of the Theorem of
 E. M. Matveev. 
\begin{thm}[Matveev]\label{thm:MM}
Let $\lambda_1$, $\lambda_2$, $\lambda_3$ be $\Q$--linearly independent logarithms of
non-zero algebraic numbers and let $b_1$, $b_2$, $b_3$ be rational integers with $b_1\not=0$.
Define $\alpha_j=\exp(\lambda_j)$ for $j=1,~ 2,~ 3$ and
\[
	\Lambda = b_1 \lambda_1 + b_2 \lambda_2 + b_3 \lambda_3 .
\]
Let $\cD$ be
the degree of the number field $\Q(\alpha_1,\alpha_2,\alpha_3)$ over $\Q$. Put 
\[
	\chi = [\R(\alpha_1,\alpha_2,\alpha_3) : \R].
\]
Let $ A_1$, $A_2$, $A_3$ 
be positive real numbers, which satisfy
\[
	A_j \ge \max\bigl\{\cD \h(\alpha_j),  |\lambda_j|, 0{.}16 \bigr\}
	\quad (1\le j\le 3) .
\]
Assume that
\[
	B \ge \max\Bigl\{1,
	\max \bigl\{ |b_j|A_j/A_1; \, 1\le j \le 3 \bigr\}  \Bigr\}.
\]
Define also
\[
	C_1 = \frac{5\times 16^5}{6\chi}\,e^3\,(7+2\chi)\left(\frac{3e}{2}\right)^\chi\!
	\Bigl( 20{.}2+\log\bigl(3^{5{.}5} \cD^2\log (e \cD )\bigr)\Bigr).
\]
Then
\[
	\log |\Lambda| > - C_1\,  \cD^2\,  A_1\, A_2 A_3\, 
	\log\, \bigl( 1{.}5\,e \cD B \log(e \cD )  \bigr) .
\]

For $\cD=2$ and $\chi=2$, this gives
\begin{equation}\label{eqn:M2}
	\log |\Lambda| > - 1{.}80741\times 10^{11} A_1\, A_2 A_3\,\log\, (13{.}80736 \,B ),
%\leqno (2)
\end{equation}
whereas, for $\cD=2$ and $\chi=1$, we get
\begin{equation}\label{eqn:M3}
	\log |\Lambda| > - 7{.}25354\times 10^{10} A_1\, A_2 A_3\,\log\, (13{.}80736 \,B ). 
%\leqno (3)
\end{equation}
\end{thm}
\begin{proof}
 	See \cite{Matveev}.
\end{proof}

\subsection{A Preliminary Bound for $p$}
It follows from Lemma~\ref{lem:kappa} that
\[
	k\Lambda = \kappa \log  \alpha  + p \log  \gamma +ir\pi
	= \kappa \log  \alpha  + p \log  \gamma + r\log(-1),
\]
which appears as a linear form of logarithms.
But a small transformation of this form leads to better estimates. Write
\[
	k\Lambda =\kappa \log (\eps_1 \alpha) + p \log (\eps_2 \gamma)+ir\pi
\]
where $\eps_1$ and $\eps_2=\pm 1$ are chosen in such a way that
\[
	\lvert \log (\eps_1 \alpha) \rvert< \pi/2
	\quad \mathrm{and} \quad
	\lvert \log (\eps_2 \gamma) \rvert < \pi/2,
\]
where we take the principal values for the logarithms, and where  $r\in \Z$ is such that 
$|\Lambda|$ is minimal (we keep the same notation $r$ as before for simplicity).

\begin{remark} 
Indeed, we can take any roots of unity in $\cL$ for $\eps_1$ 
and $\eps_2$. The only relevant case for our set of outstanding values of $D$ are
$D=25$, $100$, where $\cL=\Q(\sqrt{-1})$, whence we can realize
\[
	|\log (\eps_1 \alpha)|< \pi/4
	\quad \mathrm{and} \quad
	|\log (\eps_2 \gamma)|< \pi/4,
\]
%if $\bK=\bQ(\sqrt{-3})$ then we can realize
%$$
%|\log (\eps_1 \alpha)|< \pi/6
%\quad \text{and}\quad
%|\log (\eps_2 \gamma)|< \pi/6,
%$$
and we write
\[
\Lambda  = 2 \log  \alpha  + p \log  \gamma + r\log \zeta,
\]
where $\zeta=e^{i\pi/2}$.
%in the first case and $\zeta=e^{i\pi/3}$ in the second one.
\end{remark}

We now return to the general case.
By Lemma~\ref{lem:logLambda}  
%\[
%	\log |k \Lambda | \le - \log |x|+ \log(2{.}1\,k D_1 \sqrt{D_2} ),
%\]
%or, in terms of $y$,
\[
	\log |k\Lambda | \le - \frac{p}{2}\,\log y+ \log(2{.}2\, k D_1 \sqrt{D_2} ).
\]
%Notice that $\,h\le \frac{2}{\pi} \sqrt D_2 \log(4 D_2)$, hence
%\[
%	\log(2{.}2\, k D_1 \sqrt{D_2} ) 
%	\le \log \left(\frac{4{.}4}{\pi} \,D_1 D_2 \log(4 D_2)\right)
%	\le \log \left(\frac{4{.}4\, D}{\pi}\, \frac{20}{e} \,(4D)^{21/20}\right)
%	\le \frac{21}{20}\,\log(10\,D)
%\]
%and ({\bf H}) implies
%$$
%-\frac{n \log y}{2}+\log(2{.}2\, k c \sqrt{w} )\le
%-0{.}346\,n+ \frac{21}{20}\log(10\,C)\le -0{.}084\,n<-1{.}4,
%$$
%so that $\,|k \Lambda |<1/3\,$ and [even if $\bK=\bQ(\sqrt{-1})$ 
%or $\bK=\bQ(\sqrt{-3})$]
Our lower bound for $x,~y,~p$ implies that $\log | k \Lambda |$ is very
small and it is straightforward to deduce that 
\[
|r| \le \frac{p+1}{2}.
\]

We can write $k\Lambda$ in the form
\[
	k\Lambda = b_1 \lambda_1 + b_2 \lambda_2+ b_3\lambda_3
\]
with
\[
	b_1 =\kappa \ \text{($=1$ or 2)},\ \alpha_1 = \eps_1 \alpha, \quad
	b_2 =p,\ \alpha_2 = \eps_2 \gamma, \quad
	b_3 =r,\ \alpha_3 = -1
\]
and
$$
 \h( \alpha_1 )=   \frac{k}{\kappa}\log d, 
 \ \lambda_1= \log \alpha_1,\ \;
\h( \alpha_2 ) =  \frac{k\log y}{2}, 
\ |\lambda_2|<\pi/2, \ \;
\h( \alpha_3 ) = 0, 
\  \lambda_3=i\pi,
$$
except for the case $\cL=\Q(\sqrt{-1})$ 
%or $\bK=\bQ(\sqrt{-3})$ 
studied in the previous remark where
$\lambda_3=i\pi/2$.
% and $\lambda_3=i\pi/3$, respectively.

Applying Theorem~\ref{thm:MM}, we have $\cD=\chi=2$ and we can take  
\[
A_1 = \max\left\{\frac{2k\log d}{\kappa},\frac{\pi}{2}\right\}, \quad
A_2 = \max\left\{k\log y,\frac{\pi}{2}\right\}, \quad
A_3 = \pi   
\]
and  (using some change of numerotation in Theorem~\ref{thm:MM})
\[
B = p +1
\]
(this choice of $B$ is justified by the inequality $|r|\le (p+1)/2$ proved above), and we get
\[
	p \le C_2  k^2 \log(2 D_1) \log p.
\]
This implies
\[
	p \le  C_3 k^2 \log(2 D_1) \log \bigl(   k^2 \log(2 D_1)\bigr), 
\]
and thus
\begin{equation}\label{eqn:M4}
p \le  C_4 D_2 \log(2 D_1) \log \bigl(  D_2 \log(2D_1) \bigr),
%\leqno {(4)}
\end{equation}
where the constants could easily be explicated.

\begin{lem}\label{lem:prelimBound}
Suppose $D$ is one of the remaining values~(\ref{eqn:remaining})
and $(x,y,p)$ is a solution to~(\ref{eqn:grn})
missing from our Tables.
\begin{itemize}
\item If $D=7$ then  $p < 6{.}81 \times 10^{12}$.
\item Otherwise if $D$ is square-free then $p <1{.}448\times 10^{15}$.
\item For other values of $D$,  we have $p  <3{.}966\times 10^{14}$.
\end{itemize}
Thus in all cases $p < 1{.}5 \times 10^{15}$.
\end{lem}
\begin{proof}
This is a simple application of Matveev's Theorem~\ref{thm:MM}.
If $D=7$ it is easy to show that the $\alpha_0$ arising in 
Lemma~\ref{lem:kappa} is (upto conjugation) $(1+\sqrt{-7})/2$,
we know that $k=1$; thus $\N (\alpha_0)=2$ and 
$\Im (\log \alpha_0)=1.2094292028\ldots$   
%We assume (the equation has been solved for $y=2$ by Nagell, then $n\le 15$)
%$$
%n> 10^8  
%$$
%and then
%$$
%y \ge 22, 
%$$
%since $y$ is even, not equal to a power of 2, and splits completely in $\bQ(\sqrt{-7})$.
Then we can apply Theorem~\ref{thm:MM} with $D=2$, $\chi=2$ and
\[
 A_1 =  \pi/2 , \quad 
 A_2 =  \log y, \quad
\log A_3 =  \pi, \quad B=p+1. 
\]
After a few iterates we get the stated bound on $p$.

%We consider the equations $\,x^2+C=y^p\,$ for $\,2\le C \le 100$, with $C=c^2w$
%and $w$ square-free. In this range, all these
%equations have been solved completely, except maybe when $p$ divides the class-number, for $C$
%square-free and  ${}\not\equiv 7\pmod{8}$. Thus we limit our study at the values of $C$ in this
%range which are not square-free or  for which $C=w \equiv 7\pmod{8}$, and also where
%$d=\gcd\bigl(x-\sqrt{-C}, x+\sqrt{-C}\bigr)$ splits completely in~$\bK$.  
% \medskip
 
% The list of the remaining values of $C$ is   contained in
% $$
%\{ 7, 15_{(2)}, 23_{(3)}, 31_{(3)}, 39_{(4)}, 47_{(5)}, 
%55_{(4)}, 63, 71_{(7)}, 79_{(5)}, 87_{(6)}, 95_{(8)},
%18, 25, 28, 45_{(2)}, 60_{(2)}, 72, 92_{(3)}, 99, 100\}.
% $$
% (We have written the class-number of the field as an index between parentheses
% when it is ${}>1$.)
% For the values of this list, the class-number is always ${}\le 8$, the integer $k_0$ defined
% above is always ${}\le 8$, and $k\le 7$,  moreover 
%   $h\le 3$ when $C$ is not square-free. Besides, the gcd $d$ is always ${}\le 2$ when
% $C$ is square-free and ${}\le 10$ otherwise. More precisely, for non-squarefree $C$ 
% either the gcd is  ${}\le 10$ and $h=1$ or the gcd is ${}\le 2$, except when
% $C=45$ for which $h=2$ and the gcd is equal to~3. 
% 
% \medskip
 
In the application of Theorem~\ref{thm:MM}, we can take, for all the squarefree values of $D$,
\[
	A_1=\begin{cases}
	7\,\log2, & \text{if $k_0$ is odd,}\\
	8\,\log2, & \text{if $k_0$ is even,}
	\end{cases}
	\qquad 
	A_2 = \begin{cases}
	7 \,\log y, & \text{if $k_0$ is odd,}\\
	4\,\log y, & \text{if $k_0$ is even,}
	\end{cases}
	\qquad 
	A_3 = \pi,
\]
 so that
\[
	 A_1 A_2\le 49 \,\log2 \times \log y
\]
 and we get
\[
	 p<1{.}448\times 10^{15}.
\]
For all the remaining values of $D$, we can take
\[
	 A_1=\begin{cases}
  	\log 10, & \text{if $h=1$,}\\
 	\pi/2, & \text{if $h=2$,}\\
 	3\,\log 2, & \text{if $h=3$,}
 	\end{cases}
	 \qquad 
	A_2 = \begin{cases} 
	 \log y, & \text{if $h=1$,}\\
	 \log y, & \text{if $h=2$,}\\
	 3\,\log y, & \text{if $h=3$,}
 	\end{cases}
 	\qquad 
	A_3 = \pi,
\]
 so that
\[
	 A_1 A_2\le 9 \,\log2 \times \log y
\]
 and we get now
\[
 p <3{.}966\times 10^{14}.
\]
\end{proof}

\section{A new estimate on linear forms in three logarithms}\label{sec:3logs}

We present the type of linear forms in three logarithms that we shall study. 
We consider three non-zero  algebraic numbers $\alpha_1$, $\alpha_2$ and $\alpha_3$ 
and positive rational integers $b_1$, $b_2$, $b_3$ with $\gcd(b_1,b_2,b_3)=1$, and the linear form
\[	
	\Lambda = b_2 \log \alpha_2-b_1\log \alpha_1-b_3\log \alpha_3\not=0.
\]

We restrict our study to the following cases:

\begin{itemize}
\item {\bf the real case}: $\alpha_1$, $\alpha_2$ and $\alpha_3$ are real 
numbers ${}>1$, and the logarithms of the $\alpha_i$ are all real (and ${}>0$),
\item {\bf the complex case}: $\alpha_1$, $\alpha_2$ and $\alpha_3$ are complex numbers of 
modulus one, and the logarithms of the $\alpha_i$ are arbitrary determinations of the logarithm 
(then any of these determinations is purely imaginary).
\end{itemize}

In practice this restriction does not cause any inconvenience since 
\[
	|\Lambda| \ge \max \bigl\{|\Re(\Lambda)|, |\Im(\Lambda)|\bigr\},
\]
and so we can always reduce to the above cases.

Following \cite{B}, we use Laurent's method, and consider a suitable interpolation
determinant $\Delta$.

Without loss of generality, we may assume that
\[
	b_2 |\log \alpha_2|=  b_1|\log \alpha_1|+b_3|\log \alpha_3|\pm |\Lambda|.
\]
 
We shall choose rational positive integers $K$, $L$, $R$, $S$, $T$, with $K$, $L\ge 2$, 
we put $N=K^2L$ and we assume $RST\ge N$.
Let $i$ be an index such that $(k_i,m_i,\ell_i)$ runs trough all triples of integers with
$0\le k_i\le K-1$, $0\le m_i\le K-1$ and $0\le \ell_i\le L-1$. So each number $0$, \dots, $K-1$
occurs $KL$ times as a $k_i$, and similarly as an $m_i$, and each number $0$, \dots, $L-1$
occurs $K^2$ times as an $\ell_i$.

\smallskip

With the above definitions, let
\[
	\Delta = \det 
	\left\{ \binom{r_jb_2+s_jb_1}{k_i} \binom{t_jb_2+s_jb_3}{m_i}
	\alpha_1^{\ell_i r_j}\alpha_2^{\ell_i s_j}\alpha_3^{\ell_i t_j}
	\right\},
\]
where $r_j$, $s_j$, $t_j$ are non-negative integers less than $R$, $S$, $T$, respectively,
such that $(r_j,s_j,t_j)$ runs over $N$ distinct triples.

Put $\,\beta_1=b_1/b_2$, $\beta_3=b_3/b_2$. Let
\[
	\lambda_i = \ell_i - \frac{L-1}{2},\quad
	\eta_0 = \frac{R-1}{2}+\beta_1 \frac{S-1}{2},\quad
	\zeta_0 = \frac{T-1}{2}+\beta_3 \frac{S-1}{2},\quad
\]
and
\[
	b = (b_2\eta_0) (b_2\zeta_0)\left(\prod_{k=1}^{K-1}k!\right)^{-\frac{4}{K(K-1) }}.
\]
Following \cite{LMN}, Lemme 8, we can prove that
\begin{equation*}
	\begin{split}
	\log b 
	&\le \log \frac{(R-1)b_2+(S-1)b_1}{2}+\log \frac{(T-1)b_2+(S-1)b_3}{2} \\
	&\quad -2\log K + 3-\frac{2\log (2\pi K/e^{3/2})}{K-1}+
	\frac{2+6\pi^{-2}+\log K}{3K(K-1)}.
	\end{split}
\end{equation*}

Then, we have $\,\sum_{i=0}^{N-1}\lambda_i=0\,$ and (\cite{B}, formula (2.1))  %$%
\[
	\alpha_1^{\ell_i r_j}\alpha_2^{\ell_i s_j}\alpha_3^{\ell_i t_j}
	= \alpha_1^{\lambda_i(r_j+s_j\beta_1)}
	\alpha_3^{\lambda_i(t_j+s_j\beta_3)}(1+\theta_{ij}\Lambda'),
\]
where
\[
	\Lambda'=|\Lambda| \cdot
	\max \left\{
	\frac{LRe^{LR |\Lambda|/(2b_1)}}{2 \lvert b_1 \rvert}, %% put absolute values on b_i
	\frac{LSe^{LS |\Lambda|/(2b_2)}}{2 \lvert b_2 \rvert},
	\frac{LTe^{LT |\Lambda|/(2b_3)}}{2 \lvert b_3 \rvert}
	\right\},  %$%
\]
and where all $|\theta_{ij}|$ are ${}\le 1$.

\subsection{An upper bound for $|\Delta|$}

Put
\[
	M_1=\frac{L-1}{2}\sum_{j=1}^N r_j, \qquad
	M_2=\frac{L-1}{2}\sum_{j=1}^N s_j, \qquad
	M_3=\frac{L-1}{2}\sum_{j=1}^N t_j,  
\]
and
\[
	g=\frac{1}{4}-\frac{N}{12RST}, \quad
	G_1=\frac{NLR}{2} \,g, \quad
	G_2=\frac{NLS}{2} \,g, \quad
	G_3=\frac{NLT}{2} \,g,
\]
then, see \cite{BMS}:
 
\begin{prop}\label{prop:M1}
With  the previous notation, if  $K\ge 3$, $L\ge 5$ and $\Lambda'\le \rho^{-KL}$,
for some real number $\rho>1$, then
\begin{equation*}
	\begin{split}
	\log |\Delta| 
	&\le  \sum_{i=1}^3  M_i \log |\alpha_i|+\rho \sum_{i=1}^3 G_i |\log \alpha_i| 
	 +\log(N!)+N\log 2+ \frac{N}{2}\,(K-1) \log b \\
	& \kern 15mm  
	-\left( 
	\frac{NKL}{4}+\frac{NKL}{4(2K-1)} -\frac{NK}{3L}-\frac{N}{2} 
	\right)
	\log \rho + 0{.}0001 .
	\end{split}
\end{equation*}
\end{prop}

\subsection{A lower bound for $|\Delta|$}

Using a Liouville estimate as in \cite{LMN} Lemme~6, we get (see \cite{BMS}):

\begin{prop}\label{prop:M2}
If $\,\Delta \not= 0\,$ then
\begin{equation*}
	\begin{split}
	\log |\Delta| \ge &
	- \frac{D-1}{2}\,N\log N +\sum_{i=1}^3 ( M_i +G_i)\log |\alpha_i| \\
	 & -2D \sum_{i=1}^3  G_i \h(\alpha_i) -\frac{D-1}{2}(K-1)\,N\log b.
	\end{split}
\end{equation*}
\end{prop} 

\subsection{Synthesis}

We get (see again \cite{BMS}):
 
\begin{prop}\label{prop:M3}
With, the previous notation, if  $K\ge 3$, $L\ge 5$,
 $\rho>1$, and if $\Delta\not=0$ then
\[
	\Lambda' > \rho^{-KL}
\]
provided that
\begin{equation*}
	\begin{split}
	\left(\frac{KL}{2}+\frac{L}{4}-1   -\frac{2K}{3L} \right)  \log \rho  
	 \ge \;
	 & (D+1)\log N + gL(a_1R+a_2S+a_3T) \\
	 & +D (K-1) \, \log b - 2 \log(e/2),
	 \end{split}
\end{equation*}
where the $a_i$ are positive real numbers which satisfy
\[
	a_i\ge\rho|\log \alpha_i|- \log |\alpha_i|+2D \h(\alpha_i), \qquad
	i=1,~2,~3.
\]
\end{prop}

\begin{remark}
We notice that the statement of Proposition~\ref{prop:M3} is perfectly symmetric with respect to the
$b_i$'s or the $\alpha_i$'s, except for the choice of $b$. From now on we do not assume
that $b_1$ and $b_3$ are positive, but we still suppose that $b_2>0$ and that
\[
	b_2 |\log \alpha_2| =  |b_1\,\log \alpha_1|+ |b_3 \,\log \alpha_3|\pm |\Lambda|.
\]
\end{remark}
 %%%%%%%%%%%%%%%%%%%%%%%%%%%%%%%%%%%%%%%%%%%%%%%%%%%%%%%%%%%%%%

\subsection{A zero-lemma}

To conclude we need to find conditions under which one of our determinants $\Delta$ is non-zero,
a so-called {\it zero-lemma}. 
We use a zero-lemma due to M. Laurent \cite{L} which 
improves \cite{G} and provides an %$%
important improvement on the zero-lemma used in our previous paper \cite{BMS}:

\begin{prop}[M. Laurent]\label{prop:M4} 
 Suppose that $K$, $L$ are positive integers and that $\Sigma_1$, $\Sigma_2$ and
  $\Sigma_3$ are finite  subsets of $\bC^2\times \bC^*$ containing the origin
  and such that
\begin{equation*}
	\begin{cases}
		\card \{\lambda x_1 + \mu x_2\, : \,
		(x_1,x_2,y)\in \Sigma_1\} &>K, \quad \forall (\lambda, \mu )\not=(0,0),\\
		\card \{y\, : \,
		(x_1,x_2,y)\in \Sigma_1\} &>L,
	\end{cases}
	\leqno {\rm (i)}
\end{equation*}
and
\begin{equation*}
	\begin{cases}
		\card \{(\lambda x_1 + \mu x_2,y)\, : \,
		(x_1,x_2,y)\in \Sigma_2\} &>2KL, \quad \forall (\lambda, \mu )\not=(0,0),
		\\
		\card \{(x_1,x_2)\, : \,
		(x_1,x_2,y)\in \Sigma_2\} &>2K^2,
	\end{cases}
	\leqno {\rm (ii)}
\end{equation*}
and also that
\[
	\card  \,\Sigma_3 >6KL^2.
	\leqno {\rm (iii)}
\]
Then, the only polynomial $P\in \bC [X_1, X_2,Y]$ 
with $\deg_{X_i}P\le K$ for $i=1,~2$, and $\deg_Y P \le L$ 
which is zero on the set $\Sigma_1+\Sigma_2+\Sigma_3$, is the zero polynomial.
\end{prop} 

We now study the above conditions in detail. %% slight change
For $j=1,~2,~3$, we shall consider finite sets $\Sigma_j$ defined by
\[
	\Sigma_j = \bigl\{ (r+s\beta_1,t+s\beta_3,\alpha_1^r \alpha_2^s \alpha_3^t)
	\; :\;
	0\le r\le R_j,\; 0\le s\le S_j,\; 0\le t\le T_j\bigr\}
\]
where $R_j$, $S_j$ and $T_j$ are positive integers and where
\[
	\beta_1=\frac{b_1}{b_2}, \qquad \beta_3=\frac{b_3}{b_2}.
\] 

In practical examples, generally the following condition holds:
\begin{equation*}
	\begin{cases}
	\text{either  
	$\alpha_1$, $\alpha_2$ and $\alpha_3$ are multiplicatively independent,\ \ or} \\
	\text{two multiplicatively independent, the third a root of unity ${}\neq 1$.} 
	\end{cases}
	\leqno ({\bf M})
\end{equation*}

We also assume that
\[
	 \card \bigl \{(x_1,x_2)\, : \,(x_1,x_2,y)\in \Sigma_1\bigr\}=(R_1+1)(S_1+1)(T_1+1),
	\leqno ({\bf I}_1)
\]
and
\[
	 \card \bigl\{(x_1,x_2) \, : \,(x_1,x_2,y)\in \Sigma_2\bigr\}=(R_2+1)(S_2+1)(T_2+1). 
	\leqno ({\bf I}_2)
\]

Notice that if
\[
	(r+s\beta_1,t+s\beta_3,\alpha_1^r\alpha_2^s \alpha_3^t)=
	(r'+s'\beta_1,t'+s'\beta_3,\alpha_1^{r'}\alpha_2^{s'} \alpha_3^{t'})
\]
then, when  hypothesis ({\bf M}) holds, 
two pairs of the integers $(r,s,t)$ and $(r',s',t')$ are equal which clearly 
implies that indeed these triples are equal: for example if $\alpha_1$ and 
$\alpha_2$ are multiplicatively independent, then the
 equality $\alpha_1^r\alpha_2^s \alpha_3^t =
 \alpha_1^{r'}\alpha_2^{s'} \alpha_3^{t'}$ implies
$r=r'$ and $s=s'$ and then we conclude that $t=t'$ (use the relation $x_2=x_2'$). Hence
\[
	({\bf M}) \ \Longrightarrow \ \card\, \Sigma_j=(R_j+1) (S_j+1) (T_j+1), \quad j=1,\,2,\,3.
\]

The conditions of the zero-lemma are the following:

\noindent {\bf (i)} The first condition is divided into two subconditions
\[
	\card \bigl\{\lambda x_1 + \mu x_2\, : \,
	(x_1,x_2,y)\in \Sigma_1\bigr\} >K, \quad \forall (\lambda, \mu )\not=(0,0).
	\leqno {\rm (i.1)}
\]
This is the most technical of the above conditions, we study it in detail later.

The second subcondition is
\[
	\card \bigl\{y\, : \,
	(x_1,x_2,y)\in \Sigma_1\bigr\}  >L.
	\leqno {\rm (i.2)}
\]
We use now hypothesis ({\bf M}), and we also notice that the order
 in $\bC^*$ of a root of unity ${}\not=1$ is at least equal to 2 (since $\alpha_3\not=1$), thus
this condition is satisfied if
\begin{equation*}
	 \frac{2(R_1+1) (S_1+1)(T_1+1)}{W_1+1} >L,
	 \leqno {\rm (C.i.2)}
\end{equation*}
where $W_1$ is defined by 
\begin{equation*}
	W_1=\begin{cases}
	R_1, & \text{if $\alpha_1$ is a root of unity},
	\\
	S_1, & \text{if $\alpha_2$ is a root of unity},
	\\
	T_1, & \text{if $\alpha_3$ is a root of unity},
	\\
	1, & \text{otherwise}.
	\end{cases}
\end{equation*}
But see also the remark after (C.ii.1) below.

\noindent {\bf (ii)}  The second condition of the zero-lemma is also divided 
into two subconditions, the first being 
\[
	\card \bigl\{(\lambda x_1 + \mu x_2,y)\, : \,
	(x_1,x_2,y)\in \Sigma_2\bigr\}  >2KL, \quad \forall (\lambda, \mu )\not=(0,0).
	\leqno {\rm (ii.1)}
\]
We replace it by the stronger condition
\[
	\card \bigl\{y\, : \,
	(x_1,x_2,y)\in \Sigma_2\bigr\}  >2KL .
\]
Then, by the study of the case (i.2), we  see  that it is enough to suppose that
(when condition ({\bf M}) holds)
\begin{equation*}
 	\frac{(R_2+1) (S_2+1)(T_2+1)}{W_2+1} >KL,
	\leqno {\rm (C.ii.1)}
\end{equation*}
where $W_2$ is defined by  
\begin{equation*}
W_2=\begin{cases}
R_2, & \text{if $\alpha_1$ is a root of unity,}
\\
S_2, & \text{if $\alpha_2$ is a root of unity,}
\\
T_2, & \text{if $\alpha_3$ is a root of unity,}
\\
1, & \text{otherwise}.
\end{cases}
\end{equation*}

\begin{remark}
When (for example) $\alpha_3$ is a root of unity of order $\nu$, condition (C.ii.1)
can be replaced by
\[
 	\nu \,(R_2+1) (S_2+1)>2KL,
	\leqno {\rm (C'.ii.1)}
\]
and condition (C.i.2) can be replaced by
\[
 	\nu\,(R_1+1) (S_1+1) >L.
	\leqno {\rm (C'.i.2)}
\]
\end{remark}

The second subcondition of condition (ii) of the zero-lemma is
\[
	\card \bigl\{(x_1,x_2)\, : \, %% changed semi-colon to colon
	(x_1,x_2,y)\in \Sigma_2\bigr\}  >2K^2,
	\leqno {\rm (ii.2)}
\]
%%%%%%%%%%%%%
%\fo
%Since we assume that $\alpha_1$ and $\alpha_2$ are multiplicatively
%independent, 
%each value of $y$ determines the integers $r$ and $s$:
%
%-- if $\mu\not=0$, then
%$$
%\lambda x_1+\mu x_2 = (\lambda +\mu )r+\lambda s+\mu t
%$$ 
%and
%$$
%\card \{(x_1,x_2)\,;\,
%(x_1,x_2,y)\in \Sigma_1\} \ge (R_2+1)(S_2+1)(T_2+1).
%$$
%
%-- if  $\mu =0$, then clearly
%$$
%\card \{(x_1,x_2)\,;\,
%(x_1,x_2,y)\in \Sigma_1\} \ge  (R_2+1)(S_2+1).
%$$
%
%Thus this condition holds if
%$$
% (R_2+1) (S_2+1)>KL.
%\leqno (C.ii.2)
%$$
%\of
%%%%%%%%%%%%%
By (${\bf I}_2$) this condition is equivalent to
\[
 	(R_2+1) (S_2+1)(T_2+1)>2K^2.
	\leqno {\rm (C.ii.2)}
\]

\noindent {\bf (iii)} \ There is just one condition, namely
\[
	\card  \,\Sigma_3 >6KL^2.
\]
When ({\bf M}) holds, this is equivalent to
\[
 	(R_3+1)(S_3+1) (T_3+1)>6K^2L.
	\leqno {\rm (C.iii)}
\]

Now we have \lq translated \rq\ all the conditions of Proposition~\ref{prop:M4}, 
except the subcondition (i.1). We come back to this situation in the following Lemma which
brings some extra information to Proposition 3.1.1 of \cite{B}.

\begin{lem}\label{lem:M2} 
Let $A$, $B$ and $C$ be non-zero rational integers with
$\gcd(A,B,C)=1$ and let $D$ be an integer. Define
\[
	\Pi=\bigl\{(x,y,z)\in \bC^3 \, : \, Ax+By+Cz=D\bigr\} 
\]
and consider the set
\[
	 \Sigma = \bigl\{ (x,y,z)\in\Z^3 \, : \,  
	0\le x\le X,\; 0\le y\le Y,\; 0\le z\le Z\bigr\}, 
\]
where $X$, $Y$ and $Z$ are positive integers. Let
\[
	M = \card \,\bigl\{ (x,y,z)\in\Sigma \, : \,  Ax+By+Cz=D\bigr\}.
\]
Then
\[
	M\le \left(1+\left\lfloor \frac{X}{\alpha} \right\rfloor \right)
	\left(1+\left\lfloor \frac{Y}{|C|/\alpha} \right\rfloor \right)
	\quad \text{and} \quad
	M\le \left(1+\left\lfloor \frac{X}{\alpha} \right\rfloor \right)
	\left(1+\left\lfloor \frac{Z}{|B|/\alpha} \right\rfloor \right),
\]
where 
\[
	\alpha =\gcd(B,C).
\]

If we suppose that
\[
	M \ge  \max\bigl\{X+Y+1,\, Y+Z+1, \, Z+X+1 \bigr\}
\]
then
\begin{equation*}
	\begin{split}
	|A| & \le  \frac{(Y+1)(Z+1)}{M - \max\{Y,Z\}} \, , \quad
	|B|  \le  \frac{(X+1)(Z+1)}{M - \max\{X,Z\}} \, , \\
	|C| & \le  \frac{(X+1)(Y+1)}{M - \max\{X,Y\}}.
	\end{split}
\end{equation*}
\end{lem}
\begin{proof}
If the image (by the map $(x,y,z)\mapsto Ax+By+Cz$) of a point $(x,y,z)\in \Z^3$ belongs to the
plane $\Pi$ then 
\[
	A x \equiv D \pmod \alpha,
\]
where $A$ and $\alpha$ are coprime since $\gcd(A,B,C)=1$. This shows that the
number of such $x$ which satisfy $0\le x\le X$ is
\[
	\le 1+\left\lfloor \frac{X}{\alpha} \right\rfloor.
\] 
To simplify the notation we suppose for a while that $A$, $B$ and $C$ are positive.
Let now $x$ be fixed, with $0\le x\le X$, and such that the images of two points
$(x,y,z)$ and $(x,y',z')$ belong to $\Pi$. Then
\[
	B(y'-y)=C(z-z'),
\]
where we suppose (as we may) that $y$ is minimal (then $y'>y$). Hence there
exists $k\in \bN$ such that
\[
	y'-y = k (C/\alpha) \quad \mathrm{and}\quad
	z-z'= k (B/\alpha).
\]
It follows that, for $x$ fixed, the number of $(x,y,z)\in \Sigma$ whose
image belong to $\Pi$ is
\[
	\le 1+\left\lfloor \frac{Y}{C/\alpha} \right\rfloor .
\]
Hence
\[
	M\le \left(1+\left\lfloor \frac{X}{\alpha} \right\rfloor \right)
	\left(1+\left\lfloor \frac{Y}{C/\alpha} \right\rfloor \right),
\]
which proves the first inequality of the Lemma. The proof of the second one is
the same (looking at $z$).

\medskip

For $\xi \ge 1$ put
\[
	f(\xi) = \left( 1+\frac{X}{\xi}\right)
	\left( 1+\frac{\xi Y}{C}\right),
\]
then
\[
	M\le f(\alpha).
\]

Suppose now
\[
	M > \max\bigl\{X+1,\, Y+1, \, Z+1 \bigr\}.
\]
Put
\[
	\alpha_1=\max\{1,C/Y\}, \quad \alpha_2=\min\{C,X\}.
\]

\begin{itemize}
\item If $C>Y$ and $1\le \alpha <C/Y$ then we get $M\le X+1$, 
contradiction, thus
 \[
C>Y \ \Longrightarrow \ \alpha \ge \alpha_1 \ \text{and}\
f(\alpha_1)=2\left(1+\frac{XY}{C}\right).
\]

\item If $C>X$ and $\alpha >X$ then we get $M\le Y+1$, 
contradiction, thus
\[
	C > X \ \Longrightarrow \ \alpha \le \alpha_2 \ \mathrm{and}\
	f(\alpha_2)=2\left(1+\frac{XY}{C}\right).
\]

\item If $C\le \min\{X,Y\}$ then  $\alpha_1=1$ and  $\alpha_2=C$
and
\[
	f(\alpha_1)=(X+1)\left(1+ \frac{Y}{C}\right), \quad
	f(\alpha_2)= \left(1+ \frac{X}{C}\right)(Y+1).
\]
\end{itemize}

It is easy to check that $f''$ is positive and, from the previous study, it
follows that
\[
	M \le \max\bigl\{ f(\alpha_1),f(\alpha_2)\bigr\}.
\]
Considering the different cases $C>\max \{ X,Y\}$, $X\le C<Y$, 
$Y\le C<X$ and $C\le \min \{ X,Y\}$ we get always
%\[
%	|B|\le \frac{\nu}{\nu -1}\, \frac{(X+1)(Z+1)}{M}, \qquad
%	|C|\le \frac{\nu}{\nu -1}\, \frac{(X+1)(Y+1)}{M}.
%\]
%The proof of this is similar to that of the Lemma, but simpler. We omit the details.
\[
M \le \max\left\{
(X+1)\left(1+ \frac{Y}{ C}\right), \left(1+ \frac{X}{ C}\right)(Y+1)
\right\} =
\begin{cases} 
(X+1)\left(1+ \frac{Y}{ C}\right), & \text{if $X\ge Y$,}
\cr
\cr
\left(1+ \frac{X}{ C}\right)(Y+1), & \text{otherwise}.
 \end{cases}
\]

\smallskip

If $X\ge Y$ then
\[
 M \le (X+1)\left(1+ \frac{Y}{C}\right), 
\] 
which implies
\[
M-(X+1)\le \frac{Y(X+1)}{C}, \qquad  {\rm hence} \ \ C\le \frac{Y(X+1)}{M-(X+1)},
\] 
and the hypothesis $M\ge X+Y+1$ leads to
$$
C\le \frac{(X+1)(Y+1)}{M-X},
$$
otherwise ({\it i.e.}, if $X< Y$) we get
$$
 C \le \frac{(X+1)(Y+1)}{M-Y}.
$$
Finally, we always have
$$
|C|\le   \frac{(X+1)(Y+1)}{M-\max\{X,Y\}}.
$$

In the same way, considering now the $z$--coordinate, we get
$$
|B|\le   \frac{(X+1)(Z+1)}{M-\max\{X,Z\}}.
$$
Then, considering $y$ fixed, a similar argument gives
$$
|A|\le   \frac{(Y+1)(Z+1)}{M-\max\{X,Y\}}.
$$
\end{proof}

\medskip

\begin{cor} 
Let  $B$ and $C$ be non-zero rational integers with
$\gcd(B,C)=1$ and let $D$ be an integer. Define
the plane $\Pi$ (with $A=0$), $\Sigma$
and $M$ as in the above Lemma.
Then
\[
M\le (X+1) 
\left(1+\left\lfloor \frac{ Y}{|C|} \right\rfloor \right)
\quad \text{and} \quad
M\le (X+1)
\left(1+\left\lfloor \frac{ Z}{|B|} \right\rfloor \right).
\]
Moreover, if we suppose that
\[
M \ge \max\{ X+Y +1,X+Z+1\}
\]
then
\[
|B|\le  \frac{(X+1)(Z+1)}{M-X}, \qquad
|C|\le  \frac{(X+1)(Y+1)}{M-X}.
\]
\end{cor}
\begin{proof}
The proof is similar to that of the Lemma, but simpler. We omit the details.
\end{proof}

\medskip

\begin{lem}\label{lem:M3}
Let $R_1$, $S_1$ and $T_1$ be positive integers  
and consider the set
\[
	\tilde\Sigma_1 = \bigl\{ (x_1,x_2)=(r+s\beta_1, t+s\beta_3)\, : \;
	0\le r\le R_1,\; 0\le s\le S_1,\; 0\le t\le T_1\bigr\} 
\]
and assume  that
\[
	\card\,\tilde\Sigma_1 =(R_1+1)(S_1+1)(T_1+1).
\]
Let $(\lambda ,\mu)\in\bC^2\setminus\{(0,0)\} $ and let $c$ be a complex number.
Let $\chi$  be a positive real number.
 Then, for any $c$, the number $M$ of 
elements  $ (x_1,x_2)\in \tilde\Sigma_1$ such that 
\[
	\lambda x_1+ \mu x_2 = c
\]
satisfies
\begin{equation}\label{eqn:*}
	\;\; M \le \max\Bigl\{R_1+S_1+1,\, S_1+T_1+1, 
	\, R_1+T_1+1,\,\chi \bigl((R_1+1)( S_1+1)( T_1+1 )\bigr)^{1/2} \Bigr\} 
\end{equation}
--- except if, {\bf either} there exist two non-zero rational integers $r_1$ and $s_1$ such that
\[
	r_1b_2=s_1b_1
\]
with
\begin{equation*}
	\begin{split}
 	|r_1 | & \le \frac{(R_1+1)( T_1+1)}
	{\chi \Bigl((R_1+1)( S_1+1)( T_1+1) \Bigr)^{1/2}  -\max\{R_1,T_1\} } 
	\\ 
  	\text{and} \quad
	|s_1| & \le  \frac{(S_1+1)( T_1+1)}
	{\chi \Bigl((R_1+1)( S_1+1)( T_1+1) \Bigr)^{1/2}  -\max\{S_1,T_1\} },
	\end{split}
\end{equation*}
{\bf or} 
there exist rational integers  $r_1$, $s_1$, $t_1$ and $t_2$, with
$r_1s_1\not=0$, such that
\[
	(t_1b_1+r_1b_3)s_1=r_1b_2t_2, \qquad \gcd(r_1, t_1)=\gcd(s_1,t_2 )=1,
\]
which also satisfy
\begin{equation*}
	\begin{split}
	0  < |r_1s_1| & \le \delta \cdot  
	\frac{(R_1+1)( S_1+1) }
	{\chi \Bigl((R_1+1)( S_1+1)( T_1+1) \Bigr)^{1/2}  -\max\{R_1,S_1\} },
	\\
 	|s_1t_1| & \le \delta \cdot  
 	\frac{(S_1+1)( T_1+1) }
	{\chi \Bigl((R_1+1)( S_1+1)( T_1+1) \Bigr)^{1/2}  -\max\{S_1,T_1\} }, 
	\\ 
 	\text{and}\quad
 	|r_1t_2| & \le \delta \cdot  
 	\frac{(R_1+1)( T_1+1) }
	{\chi \Bigl((R_1+1)( S_1+1)( T_1+1) \Bigr)^{1/2}  -\max\{R_1,T_1\} },
	\\
	\end{split}
\end{equation*}
where
\[
	\delta = \gcd(r_1,s_1).
\]
If the previous upper bound~(\ref{eqn:*}) for $M$ holds then, for all $(\lambda,\mu)\in \bC^2
\setminus \{(0,0)\}$, we have
\begin{equation*}
	\begin{split}
	& \card  \bigl\{ \lambda x_1+\mu x_2 \,:\, (x_1,x_2)\in \tilde \Sigma_1\bigr\}\\
	\ge \; & 
	\frac{(R_1+1)(S_1+1)(T_1+1)}{
	\max
	\Bigl\{ R_1+S_1+1),\, S_1+T_1+1 , \, R_1+ T_1+1 ,\,\chi \bigl((R_1+1)( S_1+1)( T_1+1 )\bigr)^{1/2}
	\Bigr\}
	}.
	\end{split}
\end{equation*}
\end{lem}
\begin{proof}
Let
\[
	E_1  = \bigl\{ (r,s,t)\in \Z^3\,:\,
	0\le r\le R_1 ,\; 0\le s\le S_1 ,\; 0\le t\le T_1 \bigr\}.
\]
Recall the notation
\[
	x_1=r+\beta_1 s, \quad x_2=t+\beta_3 s, \quad \beta_1=\frac{b_1}{b_2}, 
	\quad \beta_3=\frac{b_3}{b_2}.
\]
For $(\lambda,\mu)\in \bC^2\setminus \{(0,0)\}$, we consider the cardinality 
\[
	N=\card \bigl \{ \lambda x_1+\mu x_2 \, 
	:\, (x_1,x_2)\in \tilde \Sigma_1\bigr\} .
\]
We put
\[
	M=
	\max_{c \in \bC} \card\,\bigl\{(x_1,x_2)\in \tilde\Sigma_1\,:\, \lambda x_1+\mu x_2=c\bigr\}
\]
and
\[
	\quad \Pi_c=\bigl\{(z_1,z_2)\in \bC^2\,:\, \lambda z_1+\mu z_2=c\bigr\}.
\]

We clearly have
\[
	N \ge \card \,\tilde\Sigma_1 /M,
\]
so that the last claim of the Lemma is proved and we may also suppose that~(\ref{eqn:*}) does not hold.

Consider a complex number $c$ such that the number of points 
$(x_1,x_2)\in \tilde\Sigma_1$ 
for which $\lambda x_1+\mu x_2\in\Pi_c$ is maximal (and so equal to~$M$). 

\begin{itemize}
\item If $\mu=0$: 
Apply the previous Corollary with $(x,y,z)\mapsto (r,s,t)$,
$(X,Y,Z)\mapsto (R_1,S_1,T_1)$, $(A,B,C)\mapsto (b_2/d_3,b_1/d_3,0)$, where
\[
	d_3=\gcd(b_1,b_2),
\]
and $(b_2/d_3,b_1/d_3)\mapsto (r_1,s_1)$. Then we get the wanted assertion (the
\lq either\rq\ case).

Now we assume $\mu\not=0$ and, to simplify the notation we take $\mu=1$.

\item If $\lambda=0$: Now, as above, 
apply the previous Corollary but with $(A,B,C)\mapsto (0,b_3/d_1,b_2/d_1)$, where
\[
	d_1=\gcd(b_2,b_3),
\]
and $(b_2/d_1,b_3/d_1)\mapsto (s_1,t_2)$. Then we get the asserted relation
\[
 (t_1b_1+r_1b_3)s_1=r_1b_2t_2
\]
with $r_1=1$ and $t_1=0$,
and the asserted bounds on $r_1$, $s_1$, $t_1$ and $t_2$.

\item If $\lambda b_1+b_3=0$: In this case
 $(A,B,C)\mapsto (-b_3/d_2,0,b_1/d_2)$, where
\[
	d_2=\gcd(b_1,b_3),
\]
and $(b_1/d_2,-b_3/d_2)\mapsto (r_1,t_1)$. Then we get the asserted relation
\[
 	(t_1b_1+r_1b_3)s_1=r_1b_2t_2
\]
with $s_1=1$ and $t_2=0$,
and the asserted bounds on $r_1$, $s_1$, $t_1$ and $t_2$.

\item If $\lambda \mu (\lambda b_1+b_3)\not=0$:
Since $M>S_1+1$, there exist two distinct triples $(r,s_0,t)$ and $(r',s_0,t')\in E $ such that
\[
	\lambda (r+\beta_1s_0)+(t+\beta_3s_0)= 
	\lambda (r'+\beta_1s_0)+(t'+\beta_3s_0),
\]
which gives
$\lambda (r'-r)=t-t'$,
where we suppose (as we may) that $r$ is minimal (then $r'>r$) and also that $r'-r>0$ is minimal.
Put $r_1=r'-r$ and $t_1=t-t'$, then
\[
	\lambda  =t_1/r_1.
\]
 
 Since $M>R_1+1$, there exist two distinct triples $(r_0,s,t)$ and 
$(r_0,s',t')\in E $ such that
\[
	t_1b_2r_0+(t_1b_1+r_1b_3)s+r_1b_2t= t_1b_2r_0+(t_1b_1+r_1b_3)s'+r_1b_2t',
\]
which gives now a relation of the form
\[
	(t_1b_1+r_1b_3) s_1=r_1b_2t_2,
\]
for which we may suppose that
\[
	\gcd(r_1,t_1)=\gcd(s_1,t_2)=1.
\]

Now we are ready to apply the above Corollary with 
\[
(A,B,C)\mapsto (t_1s_1/\delta,r_1t_2/\delta, r_1s_1/\delta),
\] 
where
\[
	\delta=\gcd(t_1s_1,r_1t_2, r_1s_1),
\]
and we get the conclusion, except that we have to prove that
$\delta=\gcd(r_1,s_1)$.

Suppose that $p$ is a prime divisor of $\delta$, then $p\mid r_1s_1$. 
If $p\nmid r_1$ then $p\mid s_1$ and $p\nmid t_1$, thus $p\nmid r_1t_1$:
contradiction.
If $p\nmid s_1$ then $p\mid r_1$ and $p\nmid t_1$, thus $p\nmid s_1t_2$:
contradiction. Hence, $p\mid r_1$ and $p\mid s_1$ and $p\nmid t_1t_2$. And now
it is easy to conclude that
\[
	\delta=\gcd(r_1,s_1).
\]
This ends the proof of the Lemma.
\end{itemize}
\end{proof}

\begin{remark}
Before leaving this Subsection, it is important to notice that the conclusion of the
zero-lemma, namely \lq \dots\ the only polynomial $P\in \bC [X_1, X_2,Y]$ 
with $\deg_{\bf X_i}P\le K$ for $i=1$, 2, and $\deg_Y P
\le L$ which is zero on the set $\Sigma_1+\Sigma_2+\Sigma_3$, is the zero polynomial\rq\
applied to the interpolation matrix considered above implies
that this  interpolation matrix is of maximal rank, which means that there exists a determinant
$\Delta$ as above which is nonzero.
\end{remark} 

\subsection{Statement of the main result: a lower bound for the linear form}

If we gather the results obtained in the previous subsections, we get the following theorem.

\begin{thm}\label{thm:M2}
We consider three non-zero  algebraic numbers $\alpha_1$, $\alpha_2$ 
and $\alpha_3$, all ${}\not=1$ which are either all real or all complex of modulus
one.  Moreover, we assume that
\begin{equation*}
	\begin{cases}
	\text{either  
	$\alpha_1$, $\alpha_2$ and $\alpha_3$ are multiplicatively independent,\ \ or} \\
	\text{two multiplicatively independent, the third a root of unity ${}\neq 1$.} 
	\end{cases}
	\leqno ({\bf M})
\end{equation*}
%\fo
%Moreover, we assume that
%$$
%\cases{
%\hbox{either the three numbers 
%$\alpha_1$, $\alpha_2$ and $\alpha_3$ are multiplicatively independent,}
%\cr
%\cr
%\hbox{or two of these numbers are multiplicatively independent
%and the third one is a root of unity ${}\not=1$.}
%}
%\leqno ({\bf M})
%$$
%\of
 We also consider three non-zero
rational integers $b_1$, $b_2$, $b_3$ with $\gcd(b_1,b_2,b_3)=1$, and the linear form
\[
	\Lambda = b_1\log \alpha_1+b_2 \log \alpha_2+b_3\log \alpha_3\not=0,
\]
where the logarithms of the $\alpha_i$ are arbitrary determinations of the logarithm, 
but which are all real or all purely imaginary.
Without loss of generality, we assume that
\[
	b_2 |\log \alpha_2| =  |b_1\,\log \alpha_1|+ |b_3 \,\log \alpha_3|\pm |\Lambda|.
\]
 Let $K$, $L$, $R$, $R_1$, $R_2$, $R_3$, $S$, $S_1$, $S_2$, $S_3$, $T$, $T_1$, $T_2$,
  $T_3$ be rational integers which are
all ${}\ge3$, with 
\[
	L\ge 5,\quad R>R_1+R_2+R_3,\quad S>S_1+S_2+S_3, \quad T>T_1+T_2+T_3.
\]
Let $\rho>1$ be a real number. Assume first that
\begin{equation}\label{eqn:o}
	\begin{split}
	 \left(\frac{KL}{2}+\frac{L}{4}-1  -\frac{2K}{3L} \right)  \log \rho   
  	\; \ge\; & (\cD+1)\log N + gL(a_1R+a_2S+a_3T) \\
	 & +\cD (K-1) \, \log b -2\log(e/2), \\
	\end{split}
	%\leqno{\rm (o)}
\end{equation}
where $\,N=K^2L$, 
$\,\cD=[\Q(\alpha_1,\alpha_2,\alpha_3) : \Q]\bigm/[\R(\alpha_1,\alpha_2,\alpha_3) : \R]$,
$e=\exp(1)$, %%!
\[
	g=\frac{1}{4}-\frac{N}{12RST}, \qquad 
	b = (b_2\eta_0) (b_2\zeta_0)\left(\prod_{k=1}^{K-1}k!\right)^{-\frac{4}{K(K-1) }}, 
\]
with
\[
	\eta_0 = \frac{R-1}{2}+  \frac{(S-1)b_1}{2b_2},\qquad 
	\zeta_0 = \frac{T-1}{2}+ \frac{(S-1)b_3}{2b_2}, 
\]
and
\[
	a_i \ge  \rho|\log \alpha_i|- \log |\alpha_i|+2 \cD \h(\alpha_i), \qquad
	i=1,~2,~3.
\]
{\bf If}, for some positive real number $\chi$,
\begin{enumerate}
\item[(i)] $\quad (R_1+1)(S_1+1)(T_1+1) > $
\[
	\quad  K \cdot 
	\max\Bigl\{R_1+S_1+1,\,S_1+T_1+1,\,R_1+T_1+1,\,
	\chi \bigl((R_1+1)( S_1+1)( T_1+1 )\bigr)^{1/2}\Bigr\}, 
\]
\item[(ii)] $\displaystyle \quad \card \,\bigl\{\alpha_1^r \alpha_2^s\alpha_3^t\, : 
	\, 0\le r\le R_1, \,0\le s\le S_1, \,0\le t\le T_1 \bigr\} > L$, 
\smallskip
\item[(iii)] $\displaystyle \quad (R_2+1)(S_2+1)(T_2+1) > 2K^2$, 
\smallskip
\item[(iv)] $\displaystyle \quad \card \,\bigl\{\alpha_1^r \alpha_2^s\alpha_3^t\, : 
	\, 0\le r\le R_2, \,0\le s\le S_2, \,0\le t\le T_2 \bigr\} > 2KL$, and
\smallskip
\item[(v)]  $\quad \displaystyle (R_3+1)(S_3+1)(T_3+1) > 6K^2L$, 
\end{enumerate} 
 
\noindent {\bf then either} 
\[
	 \Lambda' > \rho^{-KL},
\] 
 where
\[
	\Lambda'=|\Lambda| \cdot
	\max \left\{
	\frac{LRe^{LR |\Lambda|/(2b_1)}}{2 \lvert b_1 \rvert}, %% Cipo's correction
	\frac{LSe^{LS |\Lambda|/(2b_2)}}{2 \lvert b_2 \rvert}, %% change 
	\frac{LTe^{LT |\Lambda|/(2b_3)}}{2 \lvert b_3 \rvert} 
	\right\},
\]

\noindent {\bf or}
 at least one of the following conditions {\rm ({\bf {C1}})}, {\rm ({\bf {C2}})}, %%!
 {\rm ({\bf {C3}})} hold. %%!! 
\[
	|b_1|\le R_1\quad \text{and}\quad|b_2|\le S_1\quad \text{and}\quad |b_3|\le T_1,
	 \leqno{\bf (C1)}
\]
 
\[
 	|b_1|\le R_2\quad \text{and}\quad|b_2|\le S_2\quad \text{and}\quad |b_3|\le T_2,
 	 \leqno{\bf (C2)}
\]
 
 \smallskip
 
 \noindent
 {\rm ({\bf {C3}})} %%!
  \ \ {\bf either} there exist two non-zero rational integers $r_1$ and $s_1$ such that
\[
	r_1b_2=s_1b_1
\]
with
\begin{equation*}
	\begin{split}
 	|r_1 | & \le \frac{(R_1+1)( T_1+1)}
	{\chi \Bigl((R_1+1)( S_1+1)( T_1+1) \Bigr)^{1/2}  -\max\{R_1,T_1\} } 
	\\ 
  	\mathrm{and} \quad
	|s_1| & \le  \frac{(S_1+1)( T_1+1)}
	{\chi \Bigl((R_1+1)( S_1+1)( T_1+1) \Bigr)^{1/2}  -\max\{S_1,T_1\} },
	\end{split}
\end{equation*}
{\bf or} 
there exist rational integers  $r_1$, $s_1$, $t_1$ and $t_2$, with
$r_1s_1\not=0$, such that
\[
	(t_1b_1+r_1b_3)s_1=r_1b_2t_2, \qquad \gcd(r_1, t_1)=\gcd(s_1,t_2 )=1,
\]
which also satisfy
\begin{equation*}
	\begin{split}
	0  < |r_1s_1| & \le \delta \cdot  
	\frac{(R_1+1)( S_1+1) }
	{\chi \Bigl((R_1+1)( S_1+1)( T_1+1) \Bigr)^{1/2}  -\max\{R_1,S_1\} },
	\\
 	|s_1t_1| & \le \delta \cdot  
 	\frac{(S_1+1)( T_1+1) }
	{\chi \Bigl((R_1+1)( S_1+1)( T_1+1) \Bigr)^{1/2}  -\max\{S_1,T_1\} }, 
	\\ 
 	\mathrm{and}\quad
 	|r_1t_2| & \le \delta \cdot  
 	\frac{(R_1+1)( T_1+1) }
	{\chi \Bigl((R_1+1)( S_1+1)( T_1+1) \Bigr)^{1/2}  -\max\{R_1,T_1\} },
	\\
	\end{split}
\end{equation*}
where
\[
	\delta = \gcd(r_1,s_1).
\]
\end{thm}

\noindent{\bf Warning .} --- In the above theorem, the roles of $(\alpha_1,b_1)$ and
$(\alpha_2,b_2)$ are not completely symmetric. Even if we do not make
the hypothesis $a_1\ge a_3$ (and, of course, do not use it), 
in practice it is sometimes better to choose the numerotation such that
$a_1\ge a_3$, but one has also to deal with ({\bf C3}) which is also non-symmetrical\dots

\subsection{An estimate for linear forms in two logarithms}

We need to use linear forms in two logarithms in a very special situation
(related to condition {\bf (C3)} above) and it is difficult to find an easy-to-use
result for such a case. This is the reason why we write a suitable application of
\cite{LMN} in this Section. We apply the Corollary of Theorem~2 of \cite{M}:

\begin{prop}\label{prop:M5}
Consider the linear form in two logarithms
\[
	\Lambda=b_2\log\alpha_2-b_1\log\alpha_1,
\]
where $b_1$ and $b_2$ are positive integers. Suppose that
$\alpha_1$ and $\alpha_2$
are multiplicatively independent.
Put
\[ 
	\cD = [\Q(\alpha_1,\alpha_2) : \Q] \,/
	\, [\R(\alpha_1,\alpha_2) : \R] .
\]
Let $\, a_1$, $a_2$, $h$, $k$ be real positive numbers, and $\rho$ a real number,
$e^{3/2}\le \rho\le e^3$. Put $\lambda = \log\rho $, $\chi  =h/\lambda$ and suppose that $\chi \ge \chi _0$ for some number $\chi _0\ge 0$ and that
\begin{equation}\label{eqn:(3)}
h \ge
\max\left\{7{.}5,3\lambda,
\cD\left(\log \Bigl(\frac{b_1}{a_2}+\frac{b_2}{a_1}\Bigr)
+\log\lambda+1{.}285\right)+0.023\right\},
\end{equation}  
\begin{equation} \label{eqn:(4)}
a_i \ge 
\max\bigl\{ 4,\lambda,\rho \,|\log \alpha _i| - \log |\alpha _i|+2 \cD\h(\alpha_i) \bigr\}, 
\quad (i=1,2),\quad a_1a_2\ge 100.
\end{equation}

Put
\[
	v=4\chi  +4+1/\chi , \quad
	A=\max  \{ a_1,a_2\} .
\]
Then we have the lower bound
\[
	\log |\Lambda | \ge -  (C_0+c_1+c_2) (\lambda +h)^2 a_1a_2 ,
\]
 where
\[
	C_0 = \frac{1}{\lambda^3 }
	\left\{
	\left(2+\frac{1}{2\chi (\chi +1)} \right)
	\left(
	\frac{1}{3}+\sqrt {\frac{1}{9}+\frac{4\lambda}{3v}
	\Bigl(\frac{1}{a_1} + \frac{1}{a_2}\Bigr)
	+\frac{32\sqrt {2}(1+\chi )^{3/2}}{3v^2\sqrt {a_1a_2}}}\right)
	\right\}^2 
\]
and
\[
	c_1= \frac{\lambda (1.5\lambda +2h)}{(\lambda +h)^2 a_1a_2},\quad
	c_2= \frac{1.11 \lambda  \log \bigl( A(2\lambda +2h)^2\bigr)}{(\lambda +h)^2 a_1a_2}.
\]
\end{prop}
\begin{proof}
The only difference with Theorem~2 of \cite{M} is the definition of the term $h$.
Put
\[
	K_0:=\frac{1}{\lambda }\left(
	\frac{\sqrt {2+2\q}}{3}
	+\sqrt {\frac{2(1+\q)}{9} +
	\frac{2\lambda}{3}\left(\frac{1}{a_1}+\frac{1}{a_2}\right)
	+\frac{4\lambda \sqrt {2+\q}}{3\sqrt {a_1a_2}}}
	\,\right)^2 a_1a_2  
\]
and
\[
	f(x)=\log \frac{\bigl(1+\sqrt {x-1}\bigr) \sqrt {x}}{x-1}+ 
	\frac{\log x}{6x(x-1)} +
	\frac{3}{2}+\log \frac{3}{4}+\frac{\log \frac{x}{x-1}}{x-1}.
\]
Then the condition on $h$ in Theorem~2 of \cite{M} is
\[
	h \ge  \cD \left(\log \Bigl(\frac{b_1}{a_2}+\frac{b_2}{a_1}\Bigr)
	+\log \lambda +f(\lceil K_0\rceil)\right)+0.023.
\]
Here we can take $\chi _0=3$ and it is easy to check that our present hypotheses imply $K_0>195$. Since
$f(x)<1.285$ for $x\ge195$, we get the result.
\end{proof}

We notice that $c_1$ is a decreasing function of $\chi$, for $\chi\ge 1+\sqrt 3$ we
have
\[
	c_1 \le \frac{1}{2a_1 a_2}.
\]
We also have
\[
	c_2 \le  
	\frac{1{.}11 \lambda  \log \bigl( a_1a_2(\lambda +h)^2\bigr)}{(\lambda +h)^2 a_1a_2}
	= 
	\frac{1{.}11 \lambda }{\bigl((\lambda +h) \sqrt{a_1a_2}\bigr)^{9/5}}\,
	\frac{\log \bigl( a_1a_2(\lambda +h)^2\bigr)}{ \bigl((\lambda +h) 
	\sqrt{a_1a_2}\bigr)^{1/5}},
\]
hence
\[
	c_2 \le \frac{11{.}1 \lambda}{e\bigl((\lambda +h) \sqrt{a_1a_2}\bigr)^{9/5}} =
	\frac{11{.}1 }{e(1+\chi) (\lambda +h)^{4/5}  (a_1a_2)^{9/10}}
	<
	0{.}177\cdot (a_1a_2)^{-9/10} 
\]
(notice that the hypotheses of the above Proposition imply $\chi\ge 3$ and $\lambda +h \ge 9$).

Using these remarks and simplifying the expression of $C_0$ using $v\ge 16$ we get the
simpler estimate.

\begin{cor}\label{cor:M2} 
With the notation and hypotheses of the above proposition, we have
 the lower bound
\[
	\log |\Lambda | \ge -  (C_0'+c_1'+c_2') (\lambda +h)^2 a_1a_2 ,
\]
 where
\[
	C_0' = \frac{1}{ \lambda^3 }
	\left\{
	\left(2+\frac{1}{2\chi (\chi +1)} \right)
	\left(
	\frac{1}{3}+\sqrt {\frac{1}{9}+\frac{ \lambda}{12}
	\Bigl(\frac{1}{a_1} + \frac{1}{a_2}\Bigr)
	+\frac{\sqrt {2}}{3 \sqrt {a_1a_2}}}\right)
	\right\}^2 
\]
and
\[
	c_1' = \frac{1}{2a_1 a_2}, \qquad
	c_2' =  0{.}177\cdot (a_1a_2)^{-9/10}.
\]
\end{cor} 
%%%%%%%%%%%%%%%%%%%%%%%%%%%%%%%%%%%%%%%%%%%%%%%%%%%%%%%%%%%%%%%%%%%

\subsection{How to use Theorem~\ref{thm:M2}}

Here we assume that condition ({\bf M}) holds.
To apply the Theorem, we consider an integer $L\ge 5$ and  real parameters $m>0$,
$\rho>1$ (then one can define the $a_i$) and we put
\[
	K = \lfloor m L a_1 a_2 a_3\rfloor, \qquad
	\text{with $\ m a_1a_2 a_3\ge 2$.}
\]
To simplify the presentation, even if we do not really need these conditions, we also assume
\[
	m\ge 1, \qquad \mathrm{and} \qquad a_1, \; a_2, \; a_3\ge1.
\]

We define
\begin{equation*}
\begin{split}
R_1 &= \lfloor c_1 a_2 a_3\rfloor, \quad 
S_1 = \lfloor c_1 a_1 a_3\rfloor, \quad 
T_1 = \lfloor c_1 a_1 a_2\rfloor,
\\
R_2 &= \lfloor c_2 a_2 a_3\rfloor, \quad 
S_2 = \lfloor c_2 a_1 a_3\rfloor, \quad 
T_2 = \lfloor c_2 a_1 a_2\rfloor,
\\
R_3 &= \lfloor c_3 a_2 a_3\rfloor, \quad 
S_3 = \lfloor c_3 a_1 a_3\rfloor, \quad 
T_3 = \lfloor c_3 a_1 a_2\rfloor,
\end{split}
\end{equation*}
where the parameters $c_1$, $c_2$ and $c_3$ will be chosen so that the conditions
 (i) up to (v)  of the Theorem are satisfied.

 Clearly, condition (i) is satisfied if
\[
 	\bigl( c_1^3 (a_1 a_2 a_3)^2\bigr)^{1/2}\ge \chi m a_1 a_2 a_3 L,
	 \quad  
	 c_1^2 \cdot  a \ge  2 mL,\quad \mathrm{where}\ \;
	a =\min\{ a_1, a_2, a_3 \}.
\]
Condition (ii) is true when
$\,2c_1^2  a_1 a_2 a_3 \cdot \min\{ a_1, a_2, a_3\}\ge L$.
Thus, since we suppose $m$ and the $a_i$  all ${}\ge 1$, we can take
\[
%\boxit[2pt]
c_1 = \max \left\{ 
(\chi mL)^{2/3} , 
\left(\frac{2 mL}{a}\right)^{1/2}
\right\}.
\]
To satisfy (iii) and (iv) we can take
\[
	c_2 = \max \left\{2^{1/3}(mL)^{2/3},  \sqrt{m/a }\,L  \right\}.
\]
Finally, because of the hypothesis ({\bf M}), condition (v) holds for
\[
	c_3 = (6m^2)^{1/3}\,L.
\]

\begin{remark}
When $\alpha_1$, $\alpha_2$, $\alpha_3$ are multiplicatively independent
then it is enough to take $c_1$ and $c_3$ as above and
\[
	c_2 =  2^{1/3}(mL)^{2/3}.
\]
Then we have to verify the condition~(\ref{eqn:o}). When this inequality holds, one obtains
\[
	|\Lambda'|>\rho^{-KL},
\]
 and we get
\[
	\log |\Lambda|> - KL\log \rho -\log\bigl(\max\{R,S,T\}\cdot L\bigr), 
\]
except maybe if at least one of the conditions  ({\bf C1}), ({\bf C2}) or ({\bf C3}) holds.

Now consider the conditions ({\bf C1}), ({\bf C2}) and ({\bf C3}). For conditions ({\bf C1}) and ({\bf C2})
we have in particular
\[
	{\rm ({\bf C1})} \ \text{or} \ {\rm ({\bf C2})}\ \Longrightarrow \ 
	 b_2 \le \max\{S_1,S_2\}.
\]
Condition ({\bf C3}) will be studied for an example. Put
\[
r_1=\delta r'_1, \qquad s_1=\delta s'_1,
\]
where
\[
 \delta = \gcd(r_1, s_1). %%!
\]
We just notice, for the second alternative, 
\[
	(t_1b_1+r_1b_3)s_1=r_1b_2t_2, \qquad \gcd(r_1, t_1)=\gcd(s_1,t_2 )=1,
\]
that $r'_1\mid b_1$, say $b_1=r'_1 b_1'$, hence
\[
 	(t_1b_1'+ \delta b_3)s'_1= b_2t_2,   \quad \text{with}\ \; b_1=r_1 b_1'.
\]
If $t_2\not=0$ this shows that $s'_1\mid b_2$, say $b_2=s'_1 b_2'$, so that
\[
 	t_1b_1'+ \delta b_3 = b_2't_2,  \quad \text{with\ } \; b_1=r'_1 b_1', 
	 \ \text{and\ } b_2=s'_1 b_2'.
\]
\end{remark}
 
%\vfill \eject

\section{Completion of the Proof of Theorem~\ref{thm:main}}

Having given our new bounds for linear forms in three logarithms
we now use them to complete the proof of Theorem~\ref{thm:main}.
We have indeed shown in Lemma~\ref{lem:remaining} that if $(x,y,p)$
is a missing solution then $p > p_0$ where $p_0$ is 
given in Table~\ref{table:rational}. To complete
the proof it is enough to show that $p \leq p_0$.
In Section~\ref{sec:lin} we wrote down the linear
form in logarithms we obtain for each outstanding value
of $D$.  We will content ourselves by giving the details of this
calculation for $D=7$. The other cases are practically
identical (but with different constants, and a different number of
iterations).

%We work in the field $\cL=\Q(\sqrt{-7})$, 
%whose class-number is one.  We have the factorization 
%\[
%	x^2+7=(x+ \omega)(x- \omega)=y^p.
%\]
%Let 
%$$
%\gd = \gcd (x+ \omega, x- \omega),
%$$
%then $\gd$ divides $(2x)$ and $(2 \omega)$, so that $\gd$ divides $2$. 
%The case $\gd=1$ is the simplest and the equation is completely solved in this case\dots; 
%we consider the other case, then we may suppose that
% $\gd=(1+\omega)/2$ and we put $\alpha= (1+\omega)/2$.
% We notice that $\,\log \alpha=1{.}2094292028\ldots \times i$.
% 
%
%We have seen that
%$$
% x+ \omega = (\alpha /\bar \alpha ) \gamma ^p,
%$$
%where $\gamma$ is some non-zero integer of $\bK$.
%
%
We defined 
$$
\Lambda = \log \frac{x- \sqrt{-7}}{x+ \sqrt{-7}},
$$
and we have seen that 
$$
\log |\Lambda | \le - \frac{p}{2}\,\log y+ \log\bigl(2{.}2\, \sqrt {7} \bigr).
$$

Writing $\alpha_0=(1+\sqrt{-7})/2$ we saw that the linear form is given by
$$
\Lambda= 2 \log (\bar\alpha_0/\alpha_0)+p \log (\bar\gamma/\gamma)+i q\pi
$$
for some rational integer $q$ 
(which is not necessarily prime, but we have some lack of notation!) with $|q|<p$,
 and we get
$$
\log |\Lambda|> - KL\log \rho -\log\bigl( \max\{R,S,T\}\cdot L\bigr), 
$$
except maybe if at least one of the conditions ({\bf C1}), ({\bf C2}) or ({\bf C3}) holds.

We have already seen that
$$
 {\rm ({\bf C1})}\ {\rm or}\ {\rm ({\bf C2})} \ \Longrightarrow \ p \le \max\{S_1,S_2\}.
$$
Thus, if $p > \max\{S_1,S_2\}$ then ({\bf C1}) and ({\bf C2}) do not hold and 
then ({\bf C3}) holds, and by the study in the
previous section, either
$$
b_2=p \le   \frac{2(S_1+1)(T_1+1) }
{\chi\bigl((R_1+1) (S_1+1)(T_1+1)\bigr)^{1/2}-\max\{ S_1,T_1\}},
$$
or we obtain a relation
 $$
 t' b'+t''p+ q=0, \qquad \hbox{with \ $b'=1$ or 2,}
 $$
 (here we have used the fact that $p$ is prime, which implies $|s_1|=\delta=1$), where
 $$
  |t'|\le  \frac{(S_1+1)(T_1+1) }
  {\chi\bigl((R_1+1) (S_1+1)(T_1+1)\bigr)^{1/2}-\max\{ S_1,T_1\} }. 
 $$
 Hence,
 $$
 |t''| \le \frac{|q|}{p} + \frac{2 |t'|}{ p}< \frac{1}{ 2}+ \frac{2|t'|}{p},
 $$
which implies
$$
p <  \frac{4(S_1+1)(T_1+1) }
{\chi\bigl((R_1+1) (S_1+1)(T_1+1)\bigr)^{1/2}-\max\{ S_1,T_1\} } 
\qquad \mathrm{or}\qquad t''=0.
$$
If the second alternative holds then
$$
|q|\le 2 |t'| < 
 \frac{2(S_1+1)(T_1+1) }
 {\chi\bigl((R_1+1) (S_1+1)(T_1+1)\bigr)^{1/2}-\max\{ S_1,T_1\}  }
$$
and we can apply \cite{LMN} to the linear form in two logs
$$
\Lambda = \Log \alpha_1 + p\, \Log \alpha_2,
$$
taking $\alpha_2$ and $\,\Log \alpha_2=\log \alpha_2\,$ as before, but now
$$
\alpha_1 = \pm (\bar \alpha/\alpha)^2\qquad \text{and}\qquad 
\Log \alpha_1 = \log \alpha_1 + iq\pi.
$$
And we get for example (using Corollaire 1 and the notation of \cite{LMN})
$$
\log |\Lambda| > -31 \times\log A_1 \times\log A_2 
\times \bigl(\max\{21, \log p\}\bigr)^2.
$$

\bigskip

Now we proceed effectively to the computation of an upper bound for $p$. 
The first step is to recall that we have proved in 
Lemma~\ref{lem:prelimBound}, by applying Matveev's Theorem
(Theorem~\ref{thm:MM}),
that
\[
p< 6{.}81 \times 10^{12}. %% changed from 7 to 6.81 to match Lemma 13.4
\]

We then apply our Theorem~\ref{thm:M2} %% changed here as well, is this correct?
with the initial condition $p < 6{.}81\times 10^{12}$ and with the
lower bound
\[
y\ge 22;
\]
note that we do not yet assume our lower bound~(\ref{eqn:ylb}) obtained through
the modular approach. There are two reasons for this:
\begin{itemize}
\item The first reason is that we would like to demonstrate how 
powerful our new lower bound for linear forms in three logarithms is,
even without the help of the modular approach.
\item The second reason is that when we later make the 
assumption~(\ref{eqn:ylb}), and apply our lower bound for linear forms
in three logarithms, the reader will be
able to appreciate the saving brought by the \lq modular lower bound\rq\
 for $y$.
\end{itemize}
So for now we assume simply that $y \geq 22$ which can be deduced 
from the fact that $y$ is even, is not a power of $2$ and that $-7$ 
is a quadratic residue for every odd prime factor of $y$ (see \cite{Les}).
Applying Theorem~\ref{thm:M2} we get 
\[
p<  3{.}05\times 10^9
\]
with the choices $L=120$, $\rho=5$, $m= 106{.}2055121$, $\chi=0{.}4$ and
\[
R_1=46385, \ S_1=54196,\ T_1=37763,\  R_2=107649, \ S_2=125777, \ T_2=87639
\]
and
\[
 R_3=765790, \ S_3=894748,\  T_3=623444 
\]
 unless at least one of the conditions ({\bf C1}), ({\bf C2}), ({\bf C3}) holds.
For these values, it is clear that --- since we know that $p>10^8$ ---
conditions ({\bf C1}) and ({\bf C2}) do not hold~\footnote{
To be more precise we can take the above values for $S_1$, $S_2$ and $S_3$ 
independently of~$y$ but the $R_i$'s and $T_i$'s 
have to be increased for $y>22$,
as can be seen on the definition of the parameters given in the previous section
[$a_1$ and $a_3$ are independent of $y$ but not $a_2$]. 
Luckily, the larger $y$ is, the better 
our resulting estimate for ~$p$ will be 
and thus we can always replace $y$ by some lower bound for it.
}. %% made some slight changes in the footnote. 
We  also see that we must have $\,t''=0\,$ 
and then that
\[
	|q|<880,
\]
(which contains also the {\it exceptional case} $q=0$).
Using Corollaire 1 of \cite{LMN}, we get
\[
\log |\Lambda| > -31 \,\log A_1 \,\log A_2 \times \bigl(\max\{21, \log p\}\bigr)^2, %%!
\]
with (here)
\[
\log A_1=(|q|+2)\pi, \qquad \log A_2 = \frac{1}{2}\,\max  \{\pi, \log y\};
\]
which leads to 
\[
p < 8\times 10^7,
\]
in contradiction with our hypothesis $p>10^8$.  Thus we have proved that
\[
p <3{.}05\times 10^9.
\]
Now we iterate the same process, beginning with this new upper bound on $p$. After four iterates
(keeping the choices $L=120$ and $\rho=5$), we get
\[
p <1{.}11\times 10^9.
\]

The reader should compare this bound with the bound $p < 6{.}81 \times 10^{12}$
obtained by Matveev's Theorem.  

%Assuming now
%$$
%y>100,
%$$
%the same procedure gives
%$$
%p <8{.}25\times 10^8.
%$$
 
%It is not difficult to eliminate the range $22\le y \le 100$ (use the Baker--Davenport
%lemma), so that
%$$
%p < 8{.}25\times 10^8.
%$$

\medskip

We now assume our \lq modular lower\rq\ bound for $y$ in~(\ref{eqn:ylb}),
and then we obtain the much better bound for $p$ 
(taking $L=115$ and $\rho=5{.}4$)
\[
p <   3{.}94 \times 10^8. 
\]

\bigskip

We try to change a little the linear form.
As before we choose  $\eps=\pm 1$  such that the principal determination of $\log (\eps\bar\gamma/\gamma)$
has an absolute value ${}<\pi/2$   and we take now $\eps'$ such that the absolute value of the 
 the principal determination of $\log (\eps'\bar\alpha/\alpha)$ is minimal and thus
 ${}<\pi/2$. Then we have to distinguish two cases
$$
b_1=2, \ \alpha_1=\eps'\bar\alpha/\alpha, \quad
b_2=p, \ \alpha_2=\eps \bar\gamma/\gamma, \quad
b_3=q, \ \alpha_3=-1,
\leqno({\rm I})
$$
and
 $$
b_1=2, \ \alpha_1=\eps'\bar\alpha/\alpha, \quad
b_2=q, \ \alpha_2=-1, \quad
b_3=p, \ \alpha_3=\eps \bar\gamma/\gamma.
\leqno({\rm II})
$$

 \medskip 
 
 The study of the first case corresponds exactly to our above study, but
 with the better value
 $$
 |\log \alpha_1| =
 \min \left\{ \left| \log \frac{\bar \alpha}{\alpha }\right|,
 \pi -  \left| \log \frac{\bar \alpha}{\alpha }\right|\right\}.
 $$
 In this case we get now (with the same choices of $L$, $\rho$ and $\chi$ as above)
 $$
 p<1{.56}\times 10^8.
 $$
 
 \medskip 
 
Concerning case (II), we first notice that
$$
{\rm ({\bf C1})} \ {\rm or}\ {\rm ({\bf C2})} \ \Longrightarrow \ 
p \le \max\{T_1,T_2\},
$$
an implication which is essentially equivalent to the previous one whose conclusion was
$p \le \max\{S_1,S_2\}$. (Indeed, the present $T_i$'s play the role of the previous $S_i$'s,
and both are bounded independently of~$y$.)

\smallskip

Now we study condition ({\bf C3}). For the first alternative
$$
r_1b_2=s_1b_1,
$$
we get
 $$
 |q| < \frac{2\nu}{(\nu-1)\chi} \left(\frac{(S_1+1)(T_1+1)}{R_1+1} \right)^{1/2}
 $$
 and we can apply \cite{LMN} to the linear form in two logs
 $$
 \Lambda = \Bigl( \log(\eps'\bar \alpha/ \alpha)^2 + q\log(-1)\Bigr)
 +p\log(\eps\bar \gamma/\gamma),
 $$
 which works quite well. 
 
 Consider now the second alternative, which gives here
 $$
 2s't'+r't''q+r's'p=0.
 $$
 The cases $t'=0$ and $t''=0$ are very easy to treat, we omit the details and assume
 $t't''\not=0$. As before, dividing the above relation by $r'$, we obtain 
 $$
  s't'b'+ t''q+ s'p=0, \quad \hbox{with $\,b'=1\,$ or 2.}
 $$
 We can write
 $$
  t'b'+ t''q'+ p=0, \quad \hbox{with $\,q=s'q'$,}
 $$
 and
 $$
 t'\Lambda = q'\Bigl( \log\bigl(\alpha_1^{-2t''/b'}\bigr)  + s't'\log(-1)\Bigr)
 +p\log \bigl(\alpha_3^{t'}/\alpha_1^{2/b'}\bigr),
 $$
 which we can estimate as a linear form in two logs.
 Now we have to use Corollary~\ref{cor:M2} above. 
 
 We have the following data. We choose $L=115$, $\rho=5{.}5$ and $\chi=1$ and we get
 $$
 p<1{.}81\times 10^8, \quad \hbox{when ({\bf C3}) does not hold},
 $$
  with
 $$
 R_1=117653, \  S_1=31819, \ T_1=19991 
 $$
 and
\[ 
 	t_1= \Biggl\lceil 1.03 \sqrt{\frac{(S_1+1)(T_1+1)}{(R_1+1)}}\, \Biggr\rceil = 76, \qquad 
 t_2=\Biggl\lceil 1.03 \sqrt{\frac{(R_1+1)(T_1+1)}{(S_1+1)}}\, \Biggr\rceil = 276.
\] 
 % t1= ceil(1.01*sqrt((S1+1)*(T1+1)/(R1+1)))  
 % t2=ceil(1.01*sqrt((R1+1)*(T1+1)/(S1+1)))  

 Using  Corollary 14.8 for $\rho=8$,  we find $p<9\times 10^7$ when ({\bf C3}) holds. 
Finally, we have proved that
\[ 
 	p<1{.}81\times 10^8.
\]
 Notice that here the use of the Corollary 1 of 
\cite{LMN} produces a result which is too weak for our
 purpose, this is the reason why we have written %%! 
 here a special lower bound for linear forms in two logarithms.
 
\begin{remark}
We notice that without the modular lower bound for $y$ we were
able to show that $p < 1{.}11 \times 10^9$ whilst with this
modular lower bound we were able to improve this to 
$p < 1{.}81 \times 10^8$. Whilst it is certainly possible to reach
the former target with the methods of this paper, it would have taken
about $6$ times as long as it took to reach the latter. From this
it is a plausible guess that without the modular lower bound for $y$
the computational part for the 
entire proof for all the values of $2 \leq D \leq 100$ might have
taken  at least $1200$ days rather than $206$ days.
\end{remark}

\newpage

\section{Tables}

\begin{tabular}{||c|l||}
\hline\hline
$D$ & \ \ Solutions $\quad (\lvert x \rvert , \lvert y \rvert,n)$\\
\hline\hline
1 & $(0,1,n)$ \\
\hline
2 & $(5, 3, 3)$ \\ 
\hline
3 & \\ 
\hline
4 & $(2, 2, 3), (11, 5, 3)$ \\ 
\hline
5 & \\ 
\hline
6 & \\ 
\hline
7 & $(1, 2, 3), (181, 32, 3),(3,2,4),(5, 2, 5), (181, 8, 5),(11, 2, 7),(181,2,15)$ \\ 
\hline
8 & $(0, 2, 3)$ \\ 
\hline
9  & \\ 
\hline
10 & \\ 
\hline
11 & $(4, 3, 3), (58, 15, 3)$ \\ 
\hline
12 & $(2,2,4)$\\ 
\hline
13 & $(70, 17, 3)$ \\ 
\hline
14 & \\ 
\hline
15 & $(7, 4, 3),(1,2,4),(7,2,6)$ \\ 
\hline
16 & $(0,2,4),(4, 2, 5)$ \\
\hline
17 & $(8,3,4)$\\ 
\hline
18 & $(3, 3, 3), (15, 3, 5)$ \\ 
\hline
19 & $(18, 7, 3), (22434, 55, 5)$ \\ 
\hline
20 & $(14, 6, 3)$ \\ 
\hline
21 & \\ 
\hline
22 & \\ 
\hline
23 & $(2, 3, 3), (3, 2, 5), (45, 2, 11)$ \\ 
\hline
24 & \\ 
\hline
25 & $(10, 5, 3)$ \\ 
\hline
26 & $(1, 3, 3), (207, 35, 3) $ \\ 
\hline
27 & $(0, 3, 3)$ \\ 
\hline
28 & $(6, 4, 3), (22, 8, 3), (225, 37, 3), (2, 2, 5), (6,2,6), (10, 2, 7), (22,2,9), (362, 2, 17)$ \\ 
\hline
29 &     \\
\hline
30  & \\
\hline
31 & $(15,4,4),(1,2,5),(15,2,8)$  \\
\hline
32 & $(7,3,4),(0, 2, 5), (88, 6, 5)$ \\ 
\hline
33 & \\ 
\hline
34 & \\ 
\hline
35 & $(36, 11, 3)$ \\ 
\hline
36 & \\ 
\hline
37 & \\ 
\hline
38 &   \\
\hline
39 & $(5, 4, 3), (31, 10, 3), (103, 22, 3), (5,2,6)$ \\ 
\hline
40 & $(52, 14, 3)$ \\
\hline
41 & \\ 
\hline
42 & \\ 
\hline
43 & \\ 
\hline
44 & $(9, 5, 3)$ \\ 
\hline
45 & $(96, 21, 3), (6,3,4)$ \\ 
\hline
46 & \\ 
\hline
47 & $(13, 6, 3), (41, 12, 3), (500, 63, 3), (14, 3, 5), (9, 2, 7)$ \\ 
\hline
48 & $(4, 4, 3), (148, 28, 3),(4,2,6)$ \\ 
\hline
49 & $(524, 65, 3), (24,5,4)$ \\ 
\hline
50 &  \\
\hline\hline
\end{tabular}

\newpage

\begin{tabular}{||l|l|l||}
\hline\hline
$D$ & \ \ Solutions $\quad (\lvert x \rvert,\lvert y \rvert,n)$\\
\hline\hline
51 & \\ 
\hline
52  & \\ 
\hline
53  & $(26, 9, 3), (156, 29, 3), (26,3,6) $ \\
\hline
54  & $(17,7,3)$  \\
\hline
55  & $(3, 4, 3), (419, 56, 3), (3,2,6)$  \\
\hline
56  & $(76, 18, 3), (5,3,4)$ \\ 
\hline
57  &\\ 
\hline
58  &\\ 
\hline
59  &\\ 
\hline
60  & $(2, 4, 3), (1586, 136, 3), (14,4,4), (50354, 76, 5), (2,2,6), (14,2,8)$ \\ 
\hline
61  & $(8, 5, 3)$ \\ 
\hline
62  & \\ 
\hline
63  & $(1, 4, 3), (13537, 568, 3), (31,4,5), (1,2,6), (31,2,10)$ \\ 
\hline
64  & $(0, 4, 3), (0,2,6), (8, 2, 7)$ \\ 
\hline
65  & $(4,3,4)$\\ 
\hline
66  & \\ 
\hline
67  & $(110, 23, 3)$ \\ 
\hline
68  & \\ 
\hline
69  & \\ 
\hline
70  & \\ 
\hline
71  & $(21, 8, 3), (35,6,4), (46, 3, 7), (21,2,9)$ \\ 
\hline
72  & $(12, 6, 3), (3,3,4)$ \\ 
\hline
73  & \\ 
\hline
74  &  $(985, 99, 3), (13,3,5)$ \\
\hline
75  & \\ 
\hline
76  & $(7, 5, 3), (1015, 101, 3)$ \\ 
\hline
77  & $(2,3,4)$\\ 
\hline
78  & \\ 
\hline
79  & $(89,20,3), (7, 2, 7)$ \\ 
\hline
80  & $(1,3,4)$\\ 
\hline
81  & $(46, 13, 3), (0,3,4)$ \\ 
\hline
82  & \\ 
\hline
83  & $(140, 27, 3)$ \\ 
\hline
84  & \\ 
\hline
85  & \\ 
\hline
86  &  \\
\hline
87  & $(16, 7, 3), (13,4,4), (13,2,8)$ \\ 
\hline
88  & \\ 
\hline
89  & $(6, 5, 3)$ \\
\hline
90  & \\ 
\hline
91  & \\ 
\hline
92  & $(6, 2, 7), (90, 2, 13)$ \\ 
\hline
93  & \\ 
\hline
94  & \\ 
\hline
95  & $(11, 6, 3),(529, 6, 7)$ \\
\hline
96  & $(23,5,4)$\\ 
\hline
97  & $(48,7,4)$\\ 
\hline
98  & \\ 
\hline
99 & $(12, 3, 5)$ \\ 
\hline
100 & $(5, 5, 3), (30, 10, 3), (198, 34, 3), (55,5,5)$ \\
\hline\hline
\end{tabular}

\end{document}